\documentclass[conference]{IEEEtran}
\IEEEoverridecommandlockouts
\usepackage{amsmath,amssymb,amsthm,mathtools}
\usepackage{graphicx,booktabs,enumitem,cite,bm}
\usepackage{needspace}
\usepackage{balance}
\usepackage{color}
\usepackage[dvipsnames]{xcolor}
\usepackage[colorlinks=true,allcolors=blue]{hyperref}
\hypersetup{
pdftitle={Stability-by-Construction Approximate Dynamic Programming via Finite-Step Control Lyapunov Functions},
pdfauthor={Navid Noroozi},
pdfsubject={Integrated finite-step control, approximate dynamic programming, and robust model-based control},
pdfkeywords={finite-step control Lyapunov functions, approximate dynamic programming, constrained control}}
\usepackage{fancyhdr}
\theoremstyle{plain}
\newtheorem{theorem}{Theorem}
\newtheorem{lemma}{Lemma}
\newtheorem{proposition}{Proposition}
\newtheorem{corollary}{Corollary}
\newtheorem{definition}{Definition}
\newtheorem{assumption}{Assumption}
\newtheorem{problem}{Problem}
\newtheorem{remark}{Remark}

\newcounter{algorithm}
\newenvironment{controlalgorithm}[1]{%
  \par\Needspace{9\baselineskip}\medskip\refstepcounter{algorithm}%
  \noindent\textbf{Algorithm \thealgorithm: #1}\par\smallskip
}{\par\medskip}

\newcommand{\R}{\mathbb{R}}
\newcommand{\N}{\mathbb{N}}
\newcommand{\Rn}{\R^{n}}
\newcommand{\Rm}{\R^{m}}
\newcommand{\Rpos}{\R_{\ge 0}}
\newcommand{\Xset}{\mathcal{X}}
\newcommand{\Uset}{\mathcal{U}}

\newcommand{\Kc}{\mathcal{K}}
\newcommand{\Kinf}{\mathcal{K}_{\infty}}
\newcommand{\KL}{\mathcal{KL}}
\newcommand{\id}{\mathrm{id}}
\newcommand{\interior}{\mathrm{int}}

\newcommand{\norm}[1]{\lVert #1\rVert}
\newcommand{\Adm}{\mathcal{U}_{M}}

\newcommand{\AdmL}{\widehat{\mathcal{U}}_{M}}
\newcommand{\Bball}{\overline{\mathbb B}}
\newcommand{\Bell}{\mathcal{T}}
\newcommand{\Vc}{V_{c}^{\star}}
\newcommand{\Pistab}{\Pi_{\mathrm{stab}}}
\newcommand{\bu}{\mathbf{u}}
\newcommand{\bx}{\bar{x}}
\newcommand{\Fm}{F_{M}}
\newcommand{\Omegac}{\Omega_{c}}
\newcommand{\Fset}{\mathcal{F}}
\newcommand{\Vhat}{\widehat V}
\newcommand{\snorm}[1]{\norm{#1}_{\infty,\Omegac}}
\DeclareMathOperator*{\argmin}{arg\,min}

\allowdisplaybreaks

\begin{document}
\title{Stability-by-Construction Approximate Dynamic Programming via Finite-Step Control Lyapunov Functions}
\author{\IEEEauthorblockN{Navid Noroozi}
\thanks{N. Noroozi is with SIGNON Deutschland GmbH (DB InfraGO AG), Germany,
\texttt{navid.n.noroozi@deutschebahn.com}.}}
\maketitle

\begin{abstract}
We consider the constrained optimal control of deterministic discrete-time nonlinear
systems for which a finite-step control Lyapunov function is available. At each control
update, a sequence of $M$ inputs is chosen and applied in full. Admissibility requires a
decrease of the Lyapunov function after $M$ steps and bounds on the intermediate states.
For known dynamics, these conditions guarantee asymptotic stability for every admissible
choice of the input sequence. We minimize a discounted infinite-horizon cost over the
resulting feedback laws and bound the loss of performance in terms of the value-function
approximation error and the error in the numerical minimization. When only an approximate
model is available, we give sufficient conditions for the true trajectory to satisfy the
state constraints and for the origin to be practically asymptotically stable. We show how
to preserve these guarantees when updating the model from measurements with bounded
errors. Numerical examples include large-angle pendulum swing-up and planar trajectory
tracking in the presence of obstacles.
\end{abstract}


\section{Introduction}
\label{sec:introduction}
We develop a stability-by-construction approach to approximate dynamic programming (ADP)
as a step towards trustworthy reinforcement learning for constrained nonlinear systems.
In safety-critical applications, the use of an approximate value function to compute the
control input calls for guarantees of both closed-loop stability and constraint
satisfaction.
We therefore incorporate a finite-step control Lyapunov function (fs-CLF)~\cite{Noroozi.2020,Geiselhart.2014c} as hard constraint into the constrained optimal control problem.
In that way, we restrict the set of admissible input sequences to stabilizing control inputs, and to enforce the approximate value function to enter the objective used to select among them.

The main challenge is how to impose the Lyapunov condition while retaining as much freedom as
possible in the original constrained optimal control problem.
Requiring a given function $V$ to decrease at every sampling instant can exclude a maneuver that satisfies all state and input constraints but initially increases $V$.
Thus we require the decay only at the end of a block of $M$ steps.
An input sequence may thus produce a temporary increase of $V$, provided that the intermediate states satisfy the prescribed bounds and the terminal state satisfies the contraction condition.
This leads us to use so-called finite-step Lyapunov functions (fs-LFs)~\cite{Noroozi.2020,Geiselhart.2014c}.

The main advantage of fs-LFs over the classic LFs is the constructivity of their converse Lyapunov theorems.
While classical converse Lyapunov theorems are mainly existential, the significance of allowing several steps for Lyapunov decay can already be seen in the converse results e.g.,~\cite[Theorem~13]{Geiselhart.2014c}: 
if an asymptotic stability decay rate $\beta \in \KL$, see below for definition of class $\KL$ functions, satisfies $\beta(r,M)<r$ for a sufficiently large integer $M$, then any (scaled) norm, e.g. a quadratic function or a linearly scaled infinity norm, is a finite-step Lyapunov function.
This motivates us to integrate fs-CLF into the optimal control problem of interest, while the feasibility is ensured by the existence of the fs-CLF.

The idea of the integration of fs-CLF into the corresponding optimal control problem has been already proposed in the context of model predictive control (MPC) in~\cite{Noroozi.2020}.
In particular, the control design is formulated into a receding horizon optimization problem in which an fs-CLF is usded as a stage cost and its associated decay condition is encoded as hard constraint, cf.~\cite[Algorithm 11, 15, 19]{Noroozi.2020} for more details.
In distributed model predictive control, local finite-step conditions have also been used to account for discrepancies between predicted and realized neighbor trajectories~\cite{NorooziSharifi2026}.
Here we take the admissible input sequences of the centralized finite-step design as the decision variables of a discounted infinite-horizon control problem.

Lyapunov functions serve different purposes in learning-based control.
Chow et al. \cite{Chow.2018} develop Lyapunov-based safe policy and value iteration algorithms by imposing local constraints on a feedback law that bound the expected cumulative constraint cost in constrained Markov decision processes. 
In our work, each admissible input sequence satisfies the state and input constraints at every sample our fs-CLF framework ensures asymptotic stability and pointwise state and input constraint and a terminal fs-CLF decay constraint ensures (paractical) asymptotic stability satisfaction for deterministic nonlinear systems.
In~\cite{Choi.2020}, reinforcement learning is leveraged to compensate for model uncertainties in CLF, control barrier function, and other control-affine constraints of a quadratic program in a continuous-time setting.
This improves the performance of nominal model-based controllers through learned uncertainty compensation. 
In contrast, we propagate model-error bounds to establish sufficient conditions for constraint satisfaction and practical asymptotic stability and preserve these guarantees under data-driven model updates.

Confidence bounds on the dynamics is used in~\cite{Berkenkamp.2017} to certify a region of attraction with high probability and to guide the collection of further measurements. 
We use deterministic bounds on the model error and give a condition under which a newly fitted model preserves feasibility on a prescribed sublevel set.
Westenbroek et al. \cite{Westenbroek.2023} incorporate the one-step change of a CLF candidate into the stage cost to facilitate learning stabilizing feedback control laws with smaller discount factors. While their stability guarantees depend on the discount factor, the CLF candidate, and the suboptimality of the learned feedback control law, we impose finite-step Lyapunov contraction as a hard constraint, thereby separating stability from performance optimization.
In our setting, finite-step decay is a hard constraint on every input sequence, so an admissible control sequence preserves nominal stability even when the cost minimization is inexact. 
Li et al. \cite{Li.2026} incorporate a CLF and its decrease condition into the reward function to promote stabilizing behavior for a locomotion controller, without formally guaranteeing closed-loop stability.
We impose finite-step Lyapunov contraction as a hard constraint, ensuring nominal stability independently of value-function approximation and numerical optimization errors, provided admissibility is maintained.

Under controllability and stage-cost detectability assumptions, semiglobal practical stability of discounted infinite-horizon optimal control is investigated in~\cite{Postoyan.2017}.
They also give conditions for stability under near-optimal inputs.
 While their work concerns stability and sub-optimality analysis, our fs-CLF-based control design constraints ensure nominal stability independently of the discount factor and the accuracy of the cost minimization.
Long et al.~\cite{Long.2025} employ a fs-Layapunov based mechanism which require a weighted average of the Lyapunov values at the next $M$ closed-loop states to lie strictly below the current value away from the origin.
Their condition includes strict endpoint decrease when all weight is placed on the last step.
They augment the value function of a given feedback law with a learned residual and also train a feedback law and certificate jointly to certify a region of attraction.
Here we prescribe the fs-CLF independently of the approximate value function, thereby separating stability from performance optimization.
We further derive performance-loss estimates for constrained approximate value iteration and address bounded model uncertainty and certificate-preserving model updates.

For known dynamics, we characterize the optimal value function on each compact sublevel
set of the fs-CLF and derive performance bounds for constrained approximate value
iteration. These bounds distinguish the error in a fitted value update from the error in
minimizing over admissible input sequences. For an approximate model, we propagate a bound
on the one-step model error over the block and require each prediction ball to lie in the
prescribed domain. If such an admissible sequence exists at every state of the sublevel
set and the terminal error in the fs-CLF is smaller than the contraction margin at its
boundary, then the true trajectory satisfies the state constraints and the origin is
practically asymptotically stable. Finally, we use nested uncertainty sets to retain these
properties as new transition measurements become available.
We illustrate the effectiveness of our fs-CLF based ADP using two benchmarks: large-angle pendulum swing-up~\cite{GymnasiumBenchmark} and planar trajectory tracking in the presence of obstacles~\cite{Perille2020BARN}.

The Python implementation and the instructions for reproducing the numerical results presented in this paper are available at \href{https://github.com/navidnoroozi/fsclf4rl}{\texttt{fsclf4rl}}.

\section{Preliminaries}
\label{sec:prelim}
\subsection{Notation}
$\N=\{0,1,2,\ldots\}$, $\N^*=\N\setminus\{0\}$, $\Rpos=[0,\infty)$, and
$\norm{\cdot}$ is the Euclidean norm. A function
$\alpha:\Rpos\to\Rpos$ belongs to $\Kc$ if it is continuous, strictly increasing,
and $\alpha(0)=0$. It belongs to $\Kinf$ if additionally
$\alpha(s)\to\infty$ as $s\to\infty$. A continuous function
$\beta:\Rpos\times\Rpos\to\Rpos$ belongs to $\KL$ if
$\beta(\cdot,t)\in\Kc$ for each $t$, and $\beta(r,\cdot)$ is nonincreasing and
tends to zero for each fixed $r$. For $\alpha\in\Kinf$,
$\alpha^{-1}\in\Kinf$. We write $\rho<\id$ if $\rho(s)<s$ for every $s>0$,
and $\rho^k$ for $k$-fold composition, with $\rho^0=\id$.
Comparison function facts are used as in \cite{Kellett.2014}.
Let $\Bball(y,R)=\{z:\norm{z-y}\le R\}$ denote a closed Euclidean ball of radius $R$, centered at $y$.
A real function $W$ on a metric space is lower
semicontinuous if $W(\xi)\le\liminf_{\nu\to\infty}W(\xi_\nu)$ whenever
$\xi_\nu\to\xi$.

\subsection{System description and finite-step CLFs}
Consider the discrete-time control system
\begin{equation}
\label{eq:sys}
x(k+1)=g\big(x(k),u(k)\big),\qquad k\in\N,
\end{equation}
with $x(k)\in\Xset\subset\Rn$ and $u(k)\in\Uset\subset\Rm$.

Motivated by stabilization of system~\eqref{eq:sys} at the origin, we impose the following standing assumptions on the state an input space as well as on the dynamics $g$.

\begin{assumption}[Standing assumptions]
\label{ass:standing}
The set $\Xset$ is closed with $0\in\interior(\Xset)$, and $\Uset$ is
compact with $0\in\Uset$. 
The map $g:\Rn\times\Rm\to\Rn$ is continuous
and $g(0,0)=0$. 
There exist $\kappa_1,\kappa_2\in\Kc$ such that for all $(x,u) \in \Xset\times\Uset$
\begin{equation} \label{eq:K_bounded_g_open_loop}
\norm{g(x,u)}\le\kappa_1(\norm{x})+\kappa_2(\norm{u}) .    
\end{equation}
\end{assumption}

By $x(\cdot,\xi,\bf{u})$, we denote the solution to the controlled system~\eqref{eq:sys} under a control sequence $\bf{u} : \N \to \Uset$ with an initial condition $\xi\in\Xset$.

Given $M\in\N^*$, for each $k \in \N$, we refer to the set $\{j | j=kM+ l , l=0,\dots,M-1\}$ \emph{block $k$ of length $M$} or simply block $k$ if $M$ is clear from the context.
Fix $M\in\N^*$. For a control sequence $\bu=(u_0,\ldots,u_{M-1})\in\Uset^M$ and $\xi\in\Xset$, set
\begin{equation}
\label{eq:lift}
g^0(\xi,\bu)=\xi,\quad
g^{j+1}(\xi,\bu)=g\big(g^j(\xi,\bu),u_j\big)
\end{equation}
for $j=0,\ldots,M-1$, and define
\begin{equation}
\label{eq:m_step_block_sys}
\Fm(\xi,\bu):=g^M(\xi,\bu).
\end{equation}
We note that strictly speaking $g^j$ is only a function of $\xi$ and $(u_0,\dots,u_{j-1})$, but there is no benefit in making this precise notationaly so we stick to the simpler version.
From Assumption~\ref{ass:standing}, the maps $g^j$ are continuous in $(\xi,\bu)$.

Consider system~\eqref{eq:sys} and a map $\pi :\Xset \to \Uset^M$.
We interpret $\pi$ as a feedback evaluated at every block of length $M$.
Given an initial condition $x(0) = \xi \in\Xset$, the feedback $\pi$ determines a closed-loop trajectory as $x(j+1) = g^j(\xi,\pi(\xi))$ and at time $M$ we evaluate the feedback again and repeat the process.
In that way, the controlled system is described by
\begin{equation} \label{eq:lifted-cl}
\bx(k+1)=\Fm\big(\bx(k),\pi(\bx(k))\big),\quad \bx(0)=\xi,
\end{equation}
where $\bx(k)=x(kM)$.

We borrow the notion of admissible finite-step feedback from~\cite[Def.~3]{Noroozi.2020}.

\begin{definition}
\label{def:admff}
A map $\pi:\Xset\to\Uset^M$ is an admissible finite-step feedback of length $M$ if for every $\xi\in\Xset$ and $j=1,\ldots,M$, the following hold
\begin{enumerate}
\item[(i)] $g^j(\xi,\pi(\xi))\in\Xset$;
\item[(ii)] $\xi \mapsto g^j(\xi,\pi(\xi))$ is $\mathcal{K}$-bounded on $\Xset$, i.e., there exist functions $b_j\in\Kc$ such that
\begin{equation} \label{eq:K_bounded_g_closed_loop}
    \norm{g^j(\xi,\pi(\xi))}\le b_j(\norm{\xi}) .
\end{equation}
\end{enumerate}
\end{definition}

Condition (i) of Definition~\ref{def:admff} clarifies the terminology admissible as it ensures that trajectories of the closed-loop system obtained by applying the map $\pi$ to system~\eqref{eq:sys} stay in $\Xset$. Condition (ii) ensures that trajectories of system~\eqref{eq:sys} under the control law $\pi$ remain $\Kc$-bounded over each block.

\begin{definition}\label{def:asymp-stabilization}~\cite[Def.~4]{Noroozi.2020}
Let $\pi$ be an admissible finite-step feedback. 
Also, let $\bf{u}_\pi$ denote the generated control sequence by repeatedly applying $\pi$ at each block of length $M$.
The admissible finite-step feedback $\pi$ asymptotically stabilizes \eqref{eq:sys} in $\Xset$ if there exists $\beta\in\mathcal{KL}$ such that $\norm{x(k,\xi,\bu_\pi)}\leq\beta(\norm{\xi},k)$ for all $\xi:=x(0)\in\Xset$, $k\in\N$.
If $\beta(r,k)=C\sigma^k r$ for some $C\geq1$ and $\sigma\in[0,1)$, the stabilization is exponential.
\end{definition}

\begin{remark}
\label{rem:inputenvelope}
To check condition (ii) of Definition~\ref{def:admff}, one needs to know an admissible finite-step feedback $\pi$.
Alternatively to~\eqref{eq:K_bounded_g_closed_loop}, we directly use the $\mathcal{K}$-boundedness condition~\eqref{eq:K_bounded_g_open_loop} on $g$ may bound $\norm{u_j}\le\sigma_j(\norm{\xi})$, $\sigma_j\in\Kc$ which is more practical for numerical implementation purposes $b_0(r)=r$, $b_{j+1}(r)=\kappa_1(b_j(r))+\kappa_2(\sigma_j(r))$.
\end{remark}

Now we define the notion of finite-step control Lyapunov function (fs-CLF).
\begin{definition}
\label{def:fsclf}~\cite[Def.~6]{Noroozi.2020}
A continuous function $V:\Xset\to\Rpos$ is an fs-CLF of length $M \in\N^*$ for system~\eqref{eq:sys} if
\begin{enumerate}[label=(\roman*),leftmargin=2em]
\item there are $\alpha_1,\alpha_2\in\Kinf$ such that
\begin{equation}
\label{eq:properV}
\alpha_1(\norm{\xi})\le V(\xi)\le\alpha_2(\norm{\xi})
\quad\text{for all }\xi\in\Xset;
\end{equation}
\item there are $\rho\in\Kinf$, $\rho<\id$, and an admissible finite-step
feedback $\pi$ as in Definition~\ref{def:admff} satisfying
\begin{equation}
\label{eq:fscontraction}
V\big(\Fm(\xi,\pi(\xi))\big)\le\rho\big(V(\xi)\big)
\quad\text{for all }\xi\in\Xset.
\end{equation}
\end{enumerate}
\end{definition}

Bounds~\eqref{eq:properV} provide positive definiteness and properness
of $V$.
Inequality~\eqref{eq:fscontraction} imposes a decay condition on the state trajectories in sense of $V$ every $M$ steps.
We denote the set of the admissible finite-step feedbacks starting from each initial state $\xi \in\Xset$ as
\begin{multline}
\label{eq:admset}
\Adm(\xi):=\big\{\bu\in\Uset^M:
g^j(\xi,\bu)\in\Xset\ (j=1,\ldots,M),\\
\norm{g^j(\xi,\bu)}\le b_j(\norm{\xi})\ (j=1,\ldots,M-1),\\
V(\Fm(\xi,\bu))\le\rho(V(\xi))\big\}.
\end{multline}

The non-emptiness and compactness of $\Adm(\xi)$ is easily obtained by the existence of the fs-CLF $V$. In a similar way, the closedness of $\mathcal G$ is an immediate consequence of the closedness of $\Xset$ and compactness of $\Uset$.
Define
\begin{equation}
\label{eq:pistab}
\Pistab:=\{\pi:\Xset\to\Uset^M:
\pi(\xi)\in\Adm(\xi)\ \forall\xi\in\Xset\}.
\end{equation}
For $c>0$ and an fs-CLF $V$, let
\begin{equation}
\label{eq:omegac}
\Omegac:=\{\xi\in\Xset:V(\xi)\le c\},
\end{equation}
which is compact in $\Xset$.
Moreover, the set $\Omegac \subseteq \Xset$ is, in fact, forward invariant for system~\eqref{eq:lifted-cl}.
We present the resulting ensuring the stabilization of system~\eqref{eq:sys} at the origin by repeatedly applying a stabilizing feedback $\pi\in\Pistab$.

\begin{theorem}\cite[Prop.~8]{Noroozi.2020}
\label{thm:stab-by-constr}
Consider system~\eqref{eq:sys} satisfying Assumption~\ref{ass:standing}.
Let $V$ be an fs-CLF of length $M$ for system~\eqref{eq:sys}.
Let $\xi \in\Omegac \subseteq \Xset$ and choose any admissible finite-step feedback $\pi\in\Pistab$.
Then $\Xset$ is forward invariant and $\Omegac$ is forward invariant at each $k = l M$, $l \in \N$. Moreover, the admissible finite-step feedback $\pi\in\Pistab$ asymptotically stabilizes system~\eqref{eq:sys} at the origin relative to $\Xset$.
\end{theorem}

We close the preliminaries by recalling the following \emph{constructive} converse finite-step Lyapunov result~
\cite[Thm.~9]{Noroozi.2020}, \cite[Thm.~13]{Geiselhart.2014c}.

\begin{theorem}
\label{thm:converse}
Consider control system~\eqref{eq:sys} satisfying Assumption~\ref{ass:standing}.
Suppose that system~\eqref{eq:sys} is asymptotically stabilizable under an admissible finite-step feedback $\pi$ of length $M$.
Take any $\eta\in\Kinf$ and let $V(\xi)=\eta(\norm{\xi})$, $\theta(s)=\eta\circ\beta(\eta^{-1}(s),M)$.
Then $V$ is an fs-CLF with the decay rate $\rho(s)=(s+\theta(s))/2\in\Kinf$, $\rho<\id$ as in~\eqref{eq:fscontraction}.
\end{theorem}

Theorem~\ref{thm:converse} builds the core feasibility concept used in this work:
It ensures that if system~\eqref{eq:sys} is stabilizable in the sense of Definition~\ref{def:asymp-stabilization} with decay rate of $\beta(r,M)<r$, then any (linearly) scaled norm is fs-CLF.
This gives us the idea of choosing fs-CLF as a scaled norm of states and finding the corresponding $\pi$ and $M$, where the non-emptiness of the search is certified by Theorem~\ref{thm:converse}.
In \cite[Problem~10, Algorithm~11, and Proposition~13]{Noroozi.2020}, this fs-CLF is used as the stage cost and in a hard $M$-step contraction constraint to find such a finite-step feedback control $\pi$.
The current work extends the results in~\cite{Noroozi.2020}, where the system contains unmodeled dynamics.

\section{Constrained Cost Minimization Problem}
In this section, we define an infinite-horizon optimal control problem within our fs-CLF based framework which is numerically solved via ADP techniques in the next section.

Choose the running cost
\begin{equation}
\label{eq:stagecost}
\ell(x,u)=V(x)+r(u),
\end{equation}
where $r:\Uset\to\Rpos$ is continuous and $r(0)=0$.
This cost satisfies the lower bound $\ell(x,u)\ge\alpha_1(\norm{x})$.
Define the block cost
\begin{equation}
\label{eq:blockcost}
L(\xi,\bu):=\sum_{j=0}^{M-1}\ell(g^j(\xi,\bu),u_j) .
\end{equation}
Given a discount factor $\gamma\in(0,1)$ and a feedback law $\pi\in \Pistab$, the infinite-horizon cost is given by
\begin{equation}
\label{eq:Jpi}
J^\pi(\xi):=\sum_{k=0}^{\infty}
\gamma^k L\big(\bx(k),\pi(\bx(k))\big).
\end{equation}
We note that the discount acts \emph{per block}, while costs within a block are not discounted relative to one another.

Now we define the associated constrained infinite-horizon optimal control.
\begin{problem}
\label{prob:crl}

Consider system~\eqref{eq:sys} under Assumption~\ref{ass:standing}.
Let $V$ be an fs-CLF of length $M$ for system~\eqref{eq:sys}.
Recall the infinite-horizon cost~\eqref{eq:Jpi}.
Compute the optimal admissible finite-step feedback law $\pi^\star\in \Pistab$ which gives the following optimal value function
\begin{equation}
\label{eq:Vc}
\Vc(\xi):=\inf_{\pi\in\Pistab}J^\pi(\xi),\qquad \xi\in\Xset .
\end{equation}
\end{problem}

Let $B(\Omegac)$ be the vector space of all bounded real-valued functions $W$ on $\Omegac$, with norm $\snorm{W}:=\sup_{\xi\in\Omegac}|W(\xi)|$.
This is a complete normed space.
Let $C(\Omegac)\subset B(\Omegac)$ denote the set of the continuous functions. 
Also, let $\operatorname{LSC}_b(\Omegac)$ denote the set of the bounded lower-semicontinuous functions on $\Omegac$.
Note that boundedness alone does not imply either continuity or measurability.

Define $\bar r=\max_{u\in\Uset}r(u)$.
For every block starting in $\xi\in\Omegac$ with the control sequence $\bu \in\Adm$, we have
\begin{align}
\label{eq:costbound}
0\le L(\xi,\bu)&\le
\bar L_c:=c+\sum_{j=1}^{M-1}\alpha_2(b_j(R_c))+M\bar r,\\
\label{eq:constantupper}
0\le\Vc(\xi)\le J^\pi(\xi)&\le
\bar V_c:=\frac{\bar L_c}{1-\gamma},\qquad \pi\in\Pistab .
\end{align}
The second inequality~\eqref{eq:costbound} bounds all intermediate stage costs.
The inequality~\eqref{eq:constantupper} provides an upper bound on the optimal value function.

For $W\in B(\Omegac)$ define the constrained Bellman operator
\begin{equation}
\label{eq:bellman}
(\Bell W)(\xi):=\inf_{\bu\in\Adm(\xi)}
\{L(\xi,\bu)+\gamma W(\Fm(\xi,\bu))\}.
\end{equation}
As shown below by Theorem~\ref{thm:wellposed}(ii), the infimum in \eqref{eq:bellman} is attained for lower-semicontinuous $W$.
In particular, a continuous approximation $\Vhat$ of the optimal value function $\Vc$ defines a feedback law
\begin{equation}
\label{eq:feedback-law}
\widehat\pi(\xi)\in\argmin_{\bu\in\Adm(\xi)}
\{L(\xi,\bu)+\gamma\Vhat(\Fm(\xi,\bu))\}.
\end{equation}

\begin{remark}
The four functions $V$, $\Vc$, $\Vhat$, and $J^{\widehat\pi}$ have different roles: the fs-CLF, the optimal value function, its numerical approximation, and the actual cost of the computed feedback law, respectively. The last term in \eqref{eq:feedback-law} approximates the infinite remaining cost after the explicit block.
\end{remark}

Here we recall seminal well-posedness results on the Bellman's operator, associated optimality results from the DP literature~\cite[Prop. 2.6]{Bertsekas.1996} and its lower-semicontinuity~\cite[Thm. 3.3]{Feinberg.2014}.
\begin{theorem}
\label{thm:wellposed}
On every $\Omegac$, the following statements hold.
\begin{enumerate}[label=(\roman*),leftmargin=2em]
\item $\Bell:B(\Omegac)\to B(\Omegac)$ is monotone and satisfies
\begin{equation}
\label{eq:Tcontraction}
\snorm{\Bell W_1-\Bell W_2}\le\gamma\snorm{W_1-W_2}.
\end{equation}
\item If $W\in\operatorname{LSC}_b(\Omegac)$, the infimum in
\eqref{eq:bellman} is attained at every state, $\Bell W$ is lower
semicontinuous, and a minimizing feedback law exists.
\item The unique fixed point of $\Bell$ in $B(\Omegac)$ is the
restriction of $\Vc$ to $\Omegac$. It is lower semicontinuous.
For any $W_0\in B(\Omegac)$, the exact value iteration
$W_{i+1}=\Bell W_i$ satisfies
\begin{equation*}
\snorm{W_i-\Vc}\le\gamma^i\snorm{W_0-\Vc},
\qquad i\in\N .
\end{equation*}
\item A time-invariant feedback $\pi^\star\in\Pistab$ exists
on $\Xset$ and satisfies
\begin{align*}
\Vc(\xi)&=\min_{\bu\in\Adm(\xi)}
\{L(\xi,\bu)+\gamma\Vc(\Fm(\xi,\bu))\},\\
\pi^\star(\xi)&\in\argmin_{\bu\in\Adm(\xi)}
\{L(\xi,\bu)+\gamma\Vc(\Fm(\xi,\bu))\}.
\end{align*}
Its actual cost is $J^{\pi^\star}=\Vc$.
\end{enumerate}
\end{theorem}

Given a $\pi \in \Pistab$ on $\Omegac$, define
\begin{equation}
\label{eq:feedbackoperator}
(\Bell_\pi W)(\xi) := L(\xi,\pi(\xi))+\gamma W(\Fm(\xi,\pi(\xi))) ,
\end{equation}
which is monotone, is a $\gamma$-contraction on $B(\Omegac)$, i.e. it satisfies~\eqref{eq:Tcontraction}, and has fixed point $J^\pi$, cf. Theorem~\ref{thm:wellposed}.

\begin{theorem}
\label{thm:approx} \cite[Prop. 6.1]{Bertsekas.1996}
Let $\Vhat\in B(\Omegac)$ satisfy $\snorm{\Vhat-\Vc}\le\varepsilon$.
Suppose a feedback $\widehat\pi\in\Pistab$ exists and satisfies
\begin{equation}
\label{eq:optimizationgap}
0\le\Bell_{\widehat\pi}\Vhat-\Bell\Vhat\le\delta
\quad\text{on }\Omegac ,
\end{equation}
where $\delta\ge0$. Then it stabilizes the state from $\Omegac$ and
\begin{equation}
\label{eq:performance}
0\le J^{\widehat\pi}(\xi)-\Vc(\xi)
\le\frac{2\gamma\varepsilon+\delta}{1-\gamma},
\quad \xi\in\Omegac .
\end{equation}
For an exact minimizing feedback associated with a bounded
lower-semicontinuous $\Vhat$, one can take $\delta=0$.
\end{theorem}

\begin{remark}
Theorem~\ref{thm:approx} provides a performance estimate for the infinite-horizon cost of a suboptimal feedback $\widehat\pi$.
The smallness of $\varepsilon$ is not implied by choosing a function class or solving
a regression problem at finitely many states.
Section~\ref{sec:certification} gives sufficient ways to certify errors without knowing $\Vc$.    
\end{remark}

\section{Numerical Feedback Control and Approximate Value Iteration}
\label{sec:wp2}
In this section, we aim to numerically compute a feedback $\bu \in \Adm$ for any initial state $\xi \in \Omegac$, where the exact dynamics model $g$ as in~\eqref{eq:sys} is known.
Given fs-CLF, the contraction rate together with all the associated gains, the aim is to
compute admissible input blocks and improve the control
performance by approximating $\Vc$.

\subsection{Finding a feasible input block and improving its cost}
For a given $W:\Xset \to \Rpos$, write
\begin{equation}
\label{eq:Qobjective}
Q_W(\xi,\bu):=L(\xi,\bu)+\gamma W(\Fm(\xi,\bu)).
\end{equation}

\begin{problem}\label{prob:approximate_optimal_control}
Consider system~\eqref{eq:sys} satisfying Assumption~\ref{ass:standing}.
Let $V$ be an fs-CLF of length $M$ for system~\eqref{eq:sys}.
Given $\xi\in\Omegac$ and $\widehat V:\Xset \to \Rpos$, compute $\bu^\star \in \Adm(\xi)$ as an optimal solution of the following optimal control problem
\begin{equation}
\label{eq:feasprog}
\min_{\bu \in \Adm(\xi)} Q_{\widehat V}(\xi,\bu) .
\end{equation}
\end{problem}
Now we introduce the algorithm solving Problem~\ref{prob:approximate_optimal_control}.
\begin{remark}
In practice, if the solver additionally supplies a global lower bound
$\underline Q(\xi)$, define
\begin{equation}
\label{eq:localcontrolcertificate}
\begin{aligned}
\underline Q(\xi)&\le(\Bell\Vhat)(\xi),\\
\bar\delta(\xi)&:=Q_{\Vhat}(\xi,\bu_{\mathrm f}^\star)
                                      -\underline Q(\xi)\ge0.
\end{aligned}
\end{equation}
For a feasible returned solution, this proves
$0\le Q_{\Vhat}(\xi,\bu_{\mathrm f}^\star)
-(\Bell\Vhat)(\xi)\le\bar\delta(\xi)$.
Set $\bar\delta(\xi)=+\infty$ when no such lower bound is
available.
\end{remark}

\begin{controlalgorithm}{Feedback control computation with certified feasibility}
\label{alg:control}
\textit{Input Data}: Select model dynamics $g$ in~\eqref{eq:sys} with state constraints $\Xset$ and input constraints $\Uset$, $b_j\in\Kc$ in \eqref{eq:admset}, fs-CLF $V$ of horizon $M \in\N^*$ in Definition~\ref{def:fsclf} with the associated gains $\alpha_1,\alpha_2\in\Kinf$, $\rho\in\Kinf$, $\rho<\id$ as in~\eqref{eq:fscontraction}, initial state $\xi\in\Omegac$,
an approximate value function $\Vhat\in B(\Omegac)$, running cost $\ell$ as in~\eqref{eq:stagecost}, and discount factor $\gamma\in(0,1)$.
\begin{enumerate}[label=\arabic*.,leftmargin=1.8em]
\item Pick any initial candidate $\bu \in \Adm(\xi)$ at the current state $\xi\in\Omegac$.
\item If no verified feasible sequence is available, report failure and exit.
If a verified sequence is available, denote it by $\bu^{0}$.
Initialize $\bu_0:=\bu^0$ and $\bu_{\mathrm f}^\star:=\bu^0$.
Here $\bu_0$ denotes the initial full input block.

\item Starting from $\bu_{0}$, solve Problem~\ref{prob:approximate_optimal_control} for $\bu^\star$.
Take $\bu_{\mathrm f}^\star = \bu^\star$, if Problem~\ref{prob:approximate_optimal_control} is feasible.
\item Return $\bu_{\mathrm f}^\star$ or $(\bu_{\mathrm f}^\star,\bar\delta(\xi))$ if any low bound $\underline Q$ is available.

\end{enumerate}
\end{controlalgorithm}

Note that decreasing the finite objective $Q_{\Vhat}(\xi,\bu)$ in~\eqref{eq:Qobjective} does not, by itself, prove a decrease
in the actual infinite-horizon cost $J^\pi$ in \eqref{eq:Jpi}.
The performance guarantee given by Theorem~\ref{thm:approx} requires
the value approximation $\snorm{\Vhat-\Vc}\le\varepsilon$ and optimization-error bounds~\eqref{eq:optimizationgap}.

\subsection{Function classes and approximate dynamic programming updates}
Algorithm~\ref{alg:control} requires apriori knowledge of an approximate value function $\Vhat$.
Here we adopt the sampling-based fitted value function approach to iteratively compute the approximate value function $\Vhat$; cf. equations (1) to (3) in~\cite{Munos.2008}.
Choose a family
\begin{equation}
\label{eq:functionclass}
\Fset=\{W_\theta:\theta\in\Theta\}\subseteq C(\Omegac)
\end{equation}
of continuous approximations, where $\Theta$ is the specified
parameter set in $\R^d$ and $d$ is the number of fitted parameters.
Examples include polynomials of a chosen degree,
continuous piecewise-affine functions on a prescribed partition,
and neural networks with continuous activation functions.
Compactness of $\Omegac$ makes each such function bounded.
For convenience, assume $0\in\Fset$.

At iteration
$i\in\N$, choose $N_s\in\N^*$ sample states
$\xi^{(s)}\in\Omegac$, $s=1,\ldots,N_s$.
For each sample, compute $\bu^{(s)}\in\Adm(\xi^{(s)})$
using Algorithm~\ref{alg:control}. Its objective value is
the explicitly defined data value
\begin{equation}
\label{eq:targets}
y^{(s)}=Q_{\Vhat_i}(\xi^{(s)},\bu^{(s)}).
\end{equation}
A function in $\Fset$ is then fitted to these state--value pairs.
The fitting can be performed by a finite-dimensional regression
problem. If the input blocks are globally optimal, these are
exact values of $\Bell\Vhat_i$ at the selected states.

Define the value approximation error $e_i$, the update error
$\varepsilon_i$, and the Bellman residual $R_i$ by
\begin{align}
\label{eq:valueerror}
e_i&:=\snorm{\Vhat_i-\Vc},\\
\label{eq:fiterr}
\varepsilon_i&:=\snorm{\Vhat_{i+1}-\Bell\Vhat_i},\\
\label{eq:residual}
R_i&:=\snorm{\Vhat_i-\Bell\Vhat_i}.
\end{align}
The update error compares consecutive iterations.
The residual measures the failure of one fixed function to satisfy
the Bellman equation.

\begin{assumption}\label{ass:known_upper_bound_Vbar}
There exists a known bounded function $\bar V:\Omegac\to\Rpos$ satisfying $\Vc(\xi)\le\bar V(\xi)$ for all $\xi\in\Omegac$.
\end{assumption}

Replacing values below zero by zero and values above $\bar V(\xi)$ by $\bar V(\xi)$ is the pointwise projection
\begin{equation}
\label{eq:projection}
(\mathcal P_{\bar V}W)(\xi)
:=\min\{\bar V(\xi),\max\{0,W(\xi)\}\}.
\end{equation}
The interval $[0,\bar V(\xi)]$ is a known range for the optimal value at that state.

At iteration $i$, let $W_{i+1}\in\Fset$ denote the function fitted
to the pairs $(\xi^{(s)},y^{(s)})$ in~\eqref{eq:targets}, before projection.
Define the approximation used in the next control calculation by
\begin{equation}
\label{eq:projectediteration}
\Vhat_{i+1}:=\mathcal P_{\bar V}W_{i+1}.
\end{equation}

\begin{lemma}[Projection and error bounds]
\label{lem:clip}
Let Assumption~\ref{ass:known_upper_bound_Vbar} hold.
For every $W\in B(\Omegac)$,
\begin{equation}
\label{eq:cliptarget}
\snorm{\mathcal P_{\bar V}W-\Vc}\le\snorm{W-\Vc}.
\end{equation}
For $i\in\N$, suppose additionally that
\begin{equation}
\label{eq:projectionconditions}
\Bell\bar V\le\bar V,\qquad
0\le\Vhat_i\le\bar V\quad\text{on }\Omegac.
\end{equation}
Then $0\le\Bell\Vhat_i\le\bar V$ on $\Omegac$, and every
$W_{i+1}\in B(\Omegac)$ satisfies
\begin{equation}
\label{eq:clipupdate}
\snorm{\mathcal P_{\bar V}W_{i+1}-\Bell\Vhat_i}
\le\snorm{W_{i+1}-\Bell\Vhat_i}.
\end{equation}
If $\bar V,W\in C(\Omegac)$, then
$\mathcal P_{\bar V}W\in C(\Omegac)$.
\end{lemma}
\begin{proof}
Fix $\xi\in\Omegac$, $v\in[0,\bar V(\xi)]$, and $w\in\R$.
If $w<0$, then
\[
|\min\{\bar V(\xi),\max\{0,w\}\}-v|
=v\le v-w=|w-v|.
\]
If $0\le w\le\bar V(\xi)$, then $|\min\{\bar V(\xi),\max\{0,w\}\}-v|=|w-v|$.
If $w>\bar V(\xi)$, then
\begin{align*}
|\min\{\bar V(\xi),\max\{0,w\}\}-v|
&=\bar V(\xi)-v\\
&\le w-v=|w-v|.
\end{align*}
Thus
\begin{equation}
\label{eq:pointwiseprojection}
\begin{gathered}
|\min\{\bar V(\xi),\max\{0,w\}\}-v|\le|w-v|,\\
w\in\R,\quad v\in[0,\bar V(\xi)],\quad\xi\in\Omegac.
\end{gathered}
\end{equation}
By~\eqref{eq:constantupper} and
Assumption~\ref{ass:known_upper_bound_Vbar},
$0\le\Vc(\xi)\le\bar V(\xi)$.
Substituting $w=W(\xi)$ and $v=\Vc(\xi)$ into
\eqref{eq:pointwiseprojection} gives
\[
|(\mathcal P_{\bar V}W)(\xi)-\Vc(\xi)|
\le|W(\xi)-\Vc(\xi)|\qquad\forall\xi\in\Omegac.
\]
Taking the supremum over $\xi\in\Omegac$ proves~\eqref{eq:cliptarget}.

For $\xi\in\Omegac$ and $\bu\in\Adm(\xi)$,
equations~\eqref{eq:costbound} and~\eqref{eq:projectionconditions} imply
\begin{align*}
0&\le L(\xi,\bu)+\gamma\Vhat_i(\Fm(\xi,\bu))\\
 &\le L(\xi,\bu)+\gamma\bar V(\Fm(\xi,\bu)).
\end{align*}
Taking the infimum over $\bu\in\Adm(\xi)$ yields
\[
0\le(\Bell\Vhat_i)(\xi)
\le(\Bell\bar V)(\xi)\le\bar V(\xi).
\]
Apply~\eqref{eq:pointwiseprojection} with
$w=W_{i+1}(\xi)$ and $v=(\Bell\Vhat_i)(\xi)$ to obtain
\begin{align*}
|(\mathcal P_{\bar V}W_{i+1})(\xi)-(\Bell\Vhat_i)(\xi)| \le|W_{i+1}(\xi)-(\Bell\Vhat_i)(\xi)|.
\end{align*}
Taking the supremum proves~\eqref{eq:clipupdate}.

For $\bar V,W\in C(\Omegac)$, the identities $\max\{0,W\}=\tfrac12(W+|W|)$, $\mathcal P_{\bar V}W
=\tfrac12\big(\bar V+\max\{0,W\} - |\bar V-\max\{0,W\}|\big)$ express the projection using sums and absolute values of continuous functions. 
Hence, $\mathcal P_{\bar V}W\in C(\Omegac)$.
\end{proof}

We note that a bound $\bar V\ge\Vc$ alone does not imply $\Bell W\le\bar V$, so \eqref{eq:cliptarget} alone cannot prove~\eqref{eq:clipupdate}. 
A simple implementable choice is the constant
$\bar V\equiv\bar V_c$ from \eqref{eq:constantupper} because $(\Bell\bar V_c)(\xi)\le\bar L_c+\gamma\bar V_c=\bar V_c$.
Projection may leave the original family $\Fset$. When it is used,
the actual family is defined explicitly as
$\Fset_{\bar V}:=\{\mathcal P_{\bar V}W:W\in\Fset\}$.
For example, a projected polynomial need not be a polynomial.
Therefore, in subsequent verifications we use the projected function.

The following result collects standard approximate value-iteration error and performance bounds for the constrained Bellman operator
\eqref{eq:bellman}; see \cite{Munos.2007,Munos.2008} and \cite[Prop.~6.1]{Bertsekas.1996}.

\begin{theorem}
\label{thm:ccfvi}
Let $\Vhat_i\in B(\Omegac)$, $i\in\N$, be the actual
value-function approximations after any projection, and define
$e_i$ and $\varepsilon_i$ by \eqref{eq:valueerror}
and \eqref{eq:fiterr}, respectively.
For each $i\in\N$, let
$\pi_i:\Omegac\to\Uset^M$ be a feedback law satisfying
$\pi_i(\xi)\in\Adm(\xi)$ for every $\xi\in\Omegac$.
Suppose that there exists a finite $\delta_i\ge0$ such that
\begin{equation}
\label{eq:aviuniformgap}
0\le
(\Bell_{\pi_i}\Vhat_i)(\xi)-(\Bell\Vhat_i)(\xi)
\le\delta_i,
\qquad \forall\xi\in\Omegac.
\end{equation}
Then, for every $i\in\N$,
\begin{equation}
\label{eq:avirecursion}
e_i\le\gamma^ie_0+
\sum_{j=0}^{i-1}\gamma^{i-1-j}\varepsilon_j,
\end{equation}
Moreover, it holds that
\begin{align}
\label{eq:aviperformance}
0\le J^{\pi_i}(\xi)-\Vc(\xi)
&\le\frac{2\gamma e_i+\delta_i}{1-\gamma},
\qquad \xi\in\Omegac .
\end{align}
\end{theorem}

The value-approximation error estimate~\eqref{eq:aviuniformgap} follows from $\gamma$-contraction property~\eqref{eq:Tcontraction} and the standard approximate value-iteration error analysis; cf.~\cite[Lemma 4.1 and Thm. 5.2]{Munos.2007} and~\cite[Lemma 3]{Munos.2008}.
The feedback-performance estimate~\eqref{eq:avirecursion} follows from Theorem 4 with $\widehat V=\widehat V_i$, $\widehat\pi=\pi_i$, $\varepsilon=e_i$, $\delta=\delta_i$.
It can be also easily seen from~\eqref{eq:aviuniformgap} and~\eqref{eq:avirecursion} that if $\varepsilon_i\le\varepsilon$ and $\delta_i\le\delta$ for all $i\in\N$, where $\varepsilon,\delta\ge0$, then
\begin{align*}
\limsup_{i\to\infty}e_i&\le\frac{\varepsilon}{1-\gamma},
\\
\limsup_{i\to\infty}\snorm{J^{\pi_i}-\Vc} &\le\frac{2\gamma\varepsilon}{(1-\gamma)^2} +\frac{\delta}{1-\gamma}.
\end{align*}

For the fitted function $W_{i+1}$ in~\eqref{eq:projectediteration}, define $\varepsilon_i^{\mathrm{fit}}:=\snorm{W_{i+1}-\Bell\Vhat_i}$.
Equations~\eqref{eq:cliptarget} and~\eqref{eq:Tcontraction} give
\begin{align}
\label{eq:projectedvalueerror}
e_{i+1}
&=\snorm{\mathcal P_{\bar V}W_{i+1}-\Vc}\notag\\
&\le\snorm{W_{i+1}-\Vc}\notag\\
&\le\snorm{W_{i+1}-\Bell\Vhat_i}
       +\snorm{\Bell\Vhat_i-\Bell\Vc}\notag\\
&\le\varepsilon_i^{\mathrm{fit}}+\gamma e_i.
\end{align}
If~\eqref{eq:projectionconditions} also holds, then~\eqref{eq:clipupdate}
gives the additional update-error bound
\begin{equation}
\label{eq:projectedupdateerror}
\varepsilon_i
 =\snorm{\mathcal P_{\bar V}W_{i+1}-\Bell\Vhat_i}
 \le\varepsilon_i^{\mathrm{fit}}.
\end{equation}
Equation~\eqref{eq:projectedvalueerror} bounds $e_{i+1}$.
Equation~\eqref{eq:projectedupdateerror} bounds the distinct quantity
$\varepsilon_i$ in~\eqref{eq:fiterr}.

Theorem~\ref{thm:ccfvi} estimates $e_i$ from the errors of successive
updates. For a fixed approximation $W$, the following corollary instead
uses the single residual $\snorm{W-\Bell W}$.
It applies at any iteration with $W=\Vhat_i$ and $R=R_i$
from~\eqref{eq:residual}.

\begin{corollary}\label{cor:residual}
Let $W\in B(\Omegac)$ and $R\ge0$ satisfy
$\snorm{W-\Bell W}\le R$.
Then
\begin{equation}
\label{eq:residualvalue}
\snorm{W-\Vc}\le\frac{R}{1-\gamma}.
\end{equation}
If $\pi(\xi)\in\Adm(\xi)$ for every $\xi\in\Omegac$ and
$0\le\Bell_\pi W-\Bell W\le\delta$ on $\Omegac$, where $\delta\ge0$,
then
\begin{equation}
\label{eq:residualperformance}
\snorm{J^\pi-\Vc}\le \frac{2 R + \delta}{(1-\gamma)}.
\end{equation}
\end{corollary}

The bounds in Corollary~\ref{cor:residual} follow from
standard Bellman-residual arguments for discounted dynamic
programming; see, e.g., \cite[Thm.~5.3]{Munos.2007}.
In particular, the fixed-point identity
$\Bell\Vc=\Vc$ and the $\gamma$-contraction property
\eqref{eq:Tcontraction} yield
\eqref{eq:residualvalue}.
Similarly, since the feedback-specific Bellman operator
$\Bell_\pi$ is a $\gamma$-contraction with fixed point
$J^\pi$, the assumed optimization-gap bound implies $\snorm{W-\Bell_\pi W} \le \snorm{W-\Bell W} +\snorm{\Bell W-\Bell_\pi W}\le R+\delta$, and consequently $\snorm{J^\pi-W}\le\frac{R+\delta}{1-\gamma}$.
Combining this estimate with \eqref{eq:residualvalue}
by the triangle inequality establishes
\eqref{eq:residualperformance}.

\subsection{Uniform Error Certification and Structured Optimization}
\label{sec:certification}

Theorem~\ref{thm:ccfvi} requires bounds on the update error
$\varepsilon_i$ in~\eqref{eq:fiterr} and the optimization gap
in~\eqref{eq:aviuniformgap}.
Corollary~\ref{cor:residual} requires a bound on the Bellman residual
in~\eqref{eq:residual}.
Each of these quantities involves the infimum in~\eqref{eq:bellman}.
The following proposition replaces a bound on
$\snorm{B-\Bell A}$ by two inequalities involving the explicit
objective $Q_A$ in~\eqref{eq:Qobjective}.
Neither inequality requires knowledge of $\Vc$.

Define the compact set of admissible state--input-sequence pairs
\begin{equation}
\label{eq:Kset}
K_c:=\{(\xi,\bu):\xi\in\Omegac,\ \bu\in\Adm(\xi)\}.
\end{equation}
The following proposition extends the classical Bellman inequalities used in approximate linear programming~\cite[Eqs. (1)–(2)]{deFarias.2003} to two-sided uniform bounds on the constrained Bellman update.
\begin{proposition}
\label{prop:twosided}
Let $A,B\in B(\Omegac)$ and $\varepsilon\ge0$.
Let $\widetilde\pi:\Omegac\to\Uset^M$ satisfy
$\widetilde\pi(\xi)\in\Adm(\xi)$ for every $\xi\in\Omegac$.
Suppose
\begin{align}
\label{eq:certlower}
Q_A(\xi,\bu)+\varepsilon-B(\xi)&\ge0,
&& (\xi,\bu)\in K_c,\\
\label{eq:certupper}
B(\xi)+\varepsilon-Q_A(\xi,\widetilde\pi(\xi))&\ge0,
&& \xi\in\Omegac.
\end{align}
Then
\begin{equation}
\label{eq:certifieduniformerror}
\snorm{B-\Bell A}\le\varepsilon.
\end{equation}
In particular, the following statements hold.
\begin{enumerate}[label=(\roman*),leftmargin=2em]
\item If $A=\Vhat_i$ and $B=\Vhat_{i+1}$, then the update error
in~\eqref{eq:fiterr} satisfies $\varepsilon_i\le\varepsilon$.
\item If $A=B=W$, then
\begin{equation}
\label{eq:certifiedresiduals}
\snorm{W-\Bell W}\le\varepsilon,
\qquad
\snorm{W-\Bell_{\widetilde\pi}W}\le\varepsilon,
\end{equation}
and the infinite-horizon cost~\eqref{eq:Jpi} satisfies
\begin{equation}
\label{eq:witnessperformance}
\snorm{J^{\widetilde\pi}-\Vc}\le\frac{2\varepsilon}{1-\gamma}.
\end{equation}
\end{enumerate}
\end{proposition}
\begin{proof}
Fix $\xi\in\Omegac$.
Equation~\eqref{eq:certlower} gives $B(\xi)-\varepsilon\le Q_A(\xi,\bu)$ for all $\bu\in\Adm(\xi)$.
Taking the infimum over $\bu\in\Adm(\xi)$ and using \eqref{eq:bellman}, we get $B(\xi)-\varepsilon \le\inf_{\bu\in\Adm(\xi)}Q_A(\xi,\bu) =(\Bell A)(\xi)$.
Since $\widetilde\pi(\xi)\in\Adm(\xi)$, equation~\eqref{eq:certupper}
then gives
\begin{equation}
\label{eq:certificatechain}
\begin{aligned}
B(\xi)-\varepsilon
&\le(\Bell A)(\xi)\\
&\le Q_A(\xi,\widetilde\pi(\xi))
\le B(\xi)+\varepsilon.
\end{aligned}
\end{equation}
Therefore, $\snorm{B-\Bell A} =\sup_{\xi\in\Omegac}|B(\xi)-(\Bell A)(\xi)| \le\varepsilon$,
which proves~\eqref{eq:certifieduniformerror}.
Substituting $A=\Vhat_i$ and $B=\Vhat_{i+1}$ proves (i).

For (ii), set $A=B=W$ in~\eqref{eq:certificatechain}.
By~\eqref{eq:feedbackoperator},
\[
W(\xi)-\varepsilon
\le(\Bell W)(\xi)
\le(\Bell_{\widetilde\pi}W)(\xi)
\le W(\xi)+\varepsilon,
\]
which proves both inequalities in~\eqref{eq:certifiedresiduals}.
The fixed-point identity $\Bell\Vc=\Vc$ and
\eqref{eq:Tcontraction} imply
\begin{align*}
\snorm{W-\Vc}
&\le\snorm{W-\Bell W}+\snorm{\Bell W-\Bell\Vc}\\
&\le\varepsilon+\gamma\snorm{W-\Vc}.
\end{align*}
Similarly, $\Bell_{\widetilde\pi}J^{\widetilde\pi}=J^{\widetilde\pi}$
and the contraction property of~\eqref{eq:feedbackoperator} give
\begin{align*}
\snorm{J^{\widetilde\pi}-W}
&\le\snorm{\Bell_{\widetilde\pi}J^{\widetilde\pi}
                -\Bell_{\widetilde\pi}W}\\
&\quad+\snorm{\Bell_{\widetilde\pi}W-W}\\
&\le\gamma\snorm{J^{\widetilde\pi}-W}+\varepsilon.
\end{align*}
All norms are finite because $W\in B(\Omegac)$ and
\eqref{eq:constantupper} holds. Rearranging yields $\snorm{W-\Vc}\le\frac{\varepsilon}{1-\gamma}$ and $\snorm{J^{\widetilde\pi}-W}\le\frac{\varepsilon}{1-\gamma}$.
Consequently, we have
\[
\snorm{J^{\widetilde\pi}-\Vc}
\le\snorm{J^{\widetilde\pi}-W}+\snorm{W-\Vc}
\le\frac{2\varepsilon}{1-\gamma},
\]
which proves~\eqref{eq:witnessperformance}.
\end{proof}

The two inequalities~\eqref{eq:certlower} and~\eqref{eq:certupper}
have different requirements.
Equation~\eqref{eq:certlower} must hold for every admissible
input sequence at every state.
Equation~\eqref{eq:certupper} requires one specified feedback
$\widetilde\pi(\xi)\in\Adm(\xi)$ at each state.
This feedback need not minimize $Q_A(\xi,\cdot)$.
When $A=B=W$, a bound on $\snorm{W-\Bell_{\widetilde\pi}W}$
alone bounds $\snorm{W-J^{\widetilde\pi}}$ by the preceding proof.
The comparison with all admissible sequences in~\eqref{eq:certlower}
also bounds $\snorm{W-\Vc}$ and gives~\eqref{eq:witnessperformance}.

To compute bounds for Theorem~\ref{thm:ccfvi}, fix $i\in\N$,
set $A=\Vhat_i$ and $B=\Vhat_{i+1}$, and choose
$\widetilde\pi:\Omegac\to\Uset^M$ satisfying
$\widetilde\pi(\xi)\in\Adm(\xi)$ for every $\xi\in\Omegac$.
Define
\begin{align}
\label{eq:verifylower}
d^-&:=\sup_{(\xi,\bu)\in K_c}
             \{B(\xi)-Q_A(\xi,\bu)\},\\
\label{eq:verifyupper}
d^+&:=\sup_{\xi\in\Omegac}
             \{Q_A(\xi,\widetilde\pi(\xi))-B(\xi)\},\\
\label{eq:verifygap}
d^\pi&:=\sup_{(\xi,\bu)\in K_c}
 \{Q_A(\xi,\widetilde\pi(\xi))-Q_A(\xi,\bu)\}.
\end{align}
Upper bounds on $d^-$ and $d^+$ determine a value of $\varepsilon$
for~\eqref{eq:certlower}--\eqref{eq:certupper}.
An upper bound on $d^\pi$ determines the optimization-gap bound
in~\eqref{eq:aviuniformgap}.
The state $\xi$ is an optimization variable in
\eqref{eq:verifylower}--\eqref{eq:verifygap}, whereas it is fixed
in Algorithm~\ref{alg:control}.
The suprema are finite by~\eqref{eq:costbound} and the boundedness
of $A$ and $B$. Their attainment is not required.

\begin{proposition}
\label{prop:verification}
Fix $i\in\N$, $A=\Vhat_i$, and $B=\Vhat_{i+1}$.
Let $\widetilde\pi:\Omegac\to\Uset^M$ satisfy
$\widetilde\pi(\xi)\in\Adm(\xi)$ for every $\xi\in\Omegac$.
Let $d^-$, $d^+$, and $d^\pi$ be defined by
\eqref{eq:verifylower}--\eqref{eq:verifygap}.
\begin{enumerate}[label=(\roman*),leftmargin=2em]
\item If finite numbers $\bar d^-,\bar d^+\in\R$ satisfy
$d^-\le\bar d^-$ and $d^+\le\bar d^+$, then
\begin{equation}
\label{eq:verifiedupdate}
\varepsilon_i=\snorm{\Vhat_{i+1}-\Bell\Vhat_i}
\le\bar\varepsilon_i:=\max\{0,\bar d^-,\bar d^+\}.
\end{equation}
\item If a finite number $\bar d^\pi\in\R$ satisfies
$d^\pi\le\bar d^\pi$, then
\begin{equation}
\label{eq:verifiedgap}
\begin{gathered}
\bar\delta_i:=\max\{0,\bar d^\pi\},\\
0\le(\Bell_{\widetilde\pi}\Vhat_i)(\xi)
       -(\Bell\Vhat_i)(\xi)\le\bar\delta_i
\quad\forall\xi\in\Omegac.
\end{gathered}
\end{equation}
The bounds in (i) also give~\eqref{eq:verifiedgap} with the
alternative choice $\bar\delta_i=\max\{0,\bar d^-+\bar d^+\}$.
\item Suppose finite numbers $\bar d^-_{i,\mathrm{res}}$ and
$\bar d^+_{i,\mathrm{res}}$ satisfy
\begin{equation}
\label{eq:residualverificationbounds}
\begin{aligned}
\bar d^-_{i,\mathrm{res}}
&\ge\sup_{(\xi,\bu)\in K_c}
       \{\Vhat_i(\xi)-Q_{\Vhat_i}(\xi,\bu)\},\\
\bar d^+_{i,\mathrm{res}}
&\ge\sup_{\xi\in\Omegac}
       \{Q_{\Vhat_i}(\xi,\widetilde\pi(\xi))-\Vhat_i(\xi)\}.
\end{aligned}
\end{equation}
Then the residual $R_i$ in~\eqref{eq:residual} satisfies
\begin{equation}
\label{eq:verifiedresidual}
R_i\le\bar R_i:=
\max\{0,\bar d^-_{i,\mathrm{res}},\bar d^+_{i,\mathrm{res}}\}.
\end{equation}
\end{enumerate}
\end{proposition}
\begin{proof}
For (i), equations~\eqref{eq:verifylower}--\eqref{eq:verifyupper}
give, for every $\xi\in\Omegac$,
\begin{gather*}
B(\xi)-\bar d^-\le B(\xi)-d^-\le Q_A(\xi,\bu)
\quad\forall\bu\in\Adm(\xi),\\
Q_A(\xi,\widetilde\pi(\xi))-B(\xi)\le d^+\le\bar d^+.
\end{gather*}
Taking the infimum over $\bu\in\Adm(\xi)$ in the first
inequality and using $\widetilde\pi(\xi)\in\Adm(\xi)$ yields
\begin{equation}
\label{eq:verifiedinterval}
\begin{aligned}
B(\xi)-\bar d^-
&\le(\Bell A)(\xi)\\
&\le Q_A(\xi,\widetilde\pi(\xi))
\le B(\xi)+\bar d^+.
\end{aligned}
\end{equation}
In particular, we have $B(\xi)-(\Bell A)(\xi) \le\bar d^-$, $(\Bell A)(\xi)-B(\xi) \le\bar d^+$, $|B(\xi)-(\Bell A)(\xi)| \le\max\{0,\bar d^-,\bar d^+\}$.
Taking the supremum over $\xi\in\Omegac$ and substituting
$A=\Vhat_i$, $B=\Vhat_{i+1}$ proves~\eqref{eq:verifiedupdate}.

For (ii), fix $\xi\in\Omegac$.
Equations~\eqref{eq:bellman},~\eqref{eq:feedbackoperator}, and
\eqref{eq:verifygap} imply
\begin{align*}
0&\le(\Bell_{\widetilde\pi}A)(\xi)-(\Bell A)(\xi)\\
 &=Q_A(\xi,\widetilde\pi(\xi))
       -\inf_{\bu\in\Adm(\xi)}Q_A(\xi,\bu)\\
 &=\sup_{\bu\in\Adm(\xi)}
       \{Q_A(\xi,\widetilde\pi(\xi))-Q_A(\xi,\bu)\}\\
 &\le d^\pi\le\bar d^\pi\le\max\{0,\bar d^\pi\}.
\end{align*}
This proves~\eqref{eq:verifiedgap} with $A=\Vhat_i$.
Using~\eqref{eq:verifiedinterval} instead gives
\begin{align*}
0&\le(\Bell_{\widetilde\pi}A)(\xi)-(\Bell A)(\xi)\\
 &=[Q_A(\xi,\widetilde\pi(\xi))-B(\xi)]
       +[B(\xi)-(\Bell A)(\xi)]\\
 &\le\bar d^++\bar d^-
 \le\max\{0,\bar d^-+\bar d^+\},
\end{align*}
which proves the alternative bound in (ii).

For (iii), equation~\eqref{eq:residualverificationbounds} gives
\begin{gather*}
\Vhat_i(\xi)-\bar d^-_{i,\mathrm{res}}
\le Q_{\Vhat_i}(\xi,\bu)
\quad\forall(\xi,\bu)\in K_c,\\
Q_{\Vhat_i}(\xi,\widetilde\pi(\xi))
\le\Vhat_i(\xi)+\bar d^+_{i,\mathrm{res}}
\quad\forall\xi\in\Omegac.
\end{gather*}
Taking the infimum in the first inequality and applying
feasibility to the second yields $\Vhat_i(\xi)-\bar d^-_{i,\mathrm{res}}
\le(\Bell\Vhat_i)(\xi) \le Q_{\Vhat_i}(\xi,\widetilde\pi(\xi)) \le\Vhat_i(\xi)+\bar d^+_{i,\mathrm{res}}$ for every $\xi\in\Omegac$.
Hence, we get $R_i =\sup_{\xi\in\Omegac}|\Vhat_i(\xi)-(\Bell\Vhat_i)(\xi)| \le\max\{0,\bar d^-_{i,\mathrm{res}},\bar d^+_{i,\mathrm{res}}\} = \bar R_i$.
This proves~\eqref{eq:verifiedresidual}.
\end{proof}

Let $\pi$ denote the feedback returned by Algorithm~\ref{alg:control}
when initialized with $\widetilde\pi(\xi)$.
Suppose that $\pi(\xi)\in\Adm(\xi)$, $Q_A(\xi,\pi(\xi))\le Q_A(\xi,\widetilde\pi(\xi))$ for every $\xi\in\Omegac$.
Then the bounds~\eqref{eq:verifiedgap} and~\eqref{eq:verifyupper}
remain valid for the computed feedback as $0\le Q_A(\xi,\pi(\xi))-(\Bell A)(\xi) \le Q_A(\xi,\widetilde\pi(\xi))-(\Bell A)(\xi) \le\bar\delta_i$ and $Q_A(\xi,\pi(\xi))-B(\xi) \le Q_A(\xi,\widetilde\pi(\xi))-B(\xi) \le\bar d^+$.
An explicit candidate feedback is optional for implementation, but its full-domain feasibility is essential for this certification route.
A certified lower bound on $\Bell A$ can instead be combined with a uniform upper bound on the returned objective.
A feasible maximizer in~\eqref{eq:verifylower}--\eqref{eq:verifygap} is a \emph{lower} bound on that verification problem, not the
required upper bound. Valid global bounds are required.

\subsection{Approximate value iteration with explicit error bounds}
At each sample, Algorithm~\ref{alg:control} returns the value
$y^{(s)}$ in \eqref{eq:targets} and, if finite, the bound
$\bar\delta_i^{(s)}:=\bar\delta(\xi^{(s)})$ in
\eqref{eq:localcontrolcertificate}. Thus
\begin{equation}
\label{eq:sampletargetgap}
0\le y^{(s)}-(\Bell\Vhat_i)(\xi^{(s)})
                       \le\bar\delta_i^{(s)},\quad s=1,\ldots,N_s.
\end{equation}
For specified weights $w_s>0$, define the finite-data fitting
objective
\begin{equation}
\label{eq:generalfit}
\mathcal E_i(\theta):=\sum_{s=1}^{N_s}w_s
              [W_\theta(\xi^{(s)})-y^{(s)}]^2,\quad\theta\in\Theta.
\end{equation}

\begin{controlalgorithm}{Constrained approximate value iteration}
\label{alg:valueiteration}
\textit{Input Data:} The data of Algorithm~\ref{alg:control}, family $\Fset\subseteq C(\Omegac)$ as in~\eqref{eq:functionclass}, maximum iteration count $I_\mathrm{\max} \in \N^*$,
a tolerance $r_{\mathrm{tol}}\ge0$ for the Bellman residual
$R_i$'S, a sample count $N_s\in\N^*$,
and weights $w_s>0$ in \eqref{eq:generalfit}.
\begin{enumerate}[label=\arabic*.,leftmargin=1.8em]
\item Compute $\bar V_c$ from \eqref{eq:constantupper}. Set $\Vhat_0=0$, $\bar R_0=+\infty$,
and $i=0$.
\item Choose $\xi^{(s)}\in\Omegac$, $s=1,\ldots,N_s$.
For each $\xi^{(s)}\in\Omegac$, run Algorithm~\ref{alg:control}, with approximation
$\Vhat_i$, supplies $\bu^{(s)}\in\Adm(\xi^{(s)})$ and
the target values data~\eqref{eq:targets}.
Store $(\xi^{(s)},y^{(s)},\bar\delta_i^{(s)})$ with the
inequality \eqref{eq:sampletargetgap}. If Algorithm~\ref{alg:control} fails, stop this iteration
without forming $\Vhat_{i+1}$.
\item Minimize $\mathcal E_i(\theta)$ in \eqref{eq:generalfit}
over $\theta\in\Theta$ and select finite parameters
$\theta_{i+1}\in\Theta$. Set $W_{i+1}=W_{\theta_{i+1}}$.
Keep the projected function as the actual representation by setting $\Vhat_{i+1}=\mathcal P_{\bar V_c}W_{i+1}$ as in~\eqref{eq:projection}.
\item With $A=\Vhat_i$, $B=\Vhat_{i+1}$, compute upper
bounds for \eqref{eq:verifylower}--\eqref{eq:verifygap}, using
a feedback $\widetilde\pi_i$ for which
$\widetilde\pi_i(\xi)\in\Adm(\xi)$ has been proved for all
$\xi\in\Omegac$. Compute $\bar\varepsilon_i$ and
$\bar\delta_i$ by \eqref{eq:verifiedupdate}--\eqref{eq:verifiedgap}.
For each unavailable uniform certificate, store $+\infty$.
Also, record the sample loss $\mathcal E_i(\theta_{i+1})$.
\item Compute $\bar R_{i+1}$ by \eqref{eq:verifiedresidual}
with $A=B=\Vhat_{i+1}$, or set $\bar R_{i+1}=+\infty$
if these upper bounds are unavailable. Set $i\leftarrow i+1$.
Stop if $\bar R_i\le r_{\mathrm{tol}}$ or if $i = I_\mathrm{\max}$. 
Otherwise return to step 2.
\item Return $\Vhat_i$, the bounds
$\{\bar\varepsilon_j,\bar\delta_j\}_{j<i}$, and $\bar R_i$.
The implemented feedback at a measured state is the input
block returned by Algorithm~\ref{alg:control} using $\Vhat_i$.
If a uniform gap $\bar\delta_i$ for this final feedback is
needed, certify \eqref{eq:verifygap} with $A=\Vhat_i$.
\end{enumerate}
\end{controlalgorithm}
Finite values in \eqref{eq:verifiedupdate}--\eqref{eq:verifiedgap}
can be substituted into \eqref{eq:avirecursion} and
\eqref{eq:aviperformance}. The stopping test proves
$R_i\le r_{\mathrm{tol}}$ without knowing $\Vc$.
Moreover, \eqref{eq:residualvalue}--\eqref{eq:residualperformance}
 give the associated value and performance bounds.
If a certificate is $+\infty$, no finite performance bound
is asserted from it. Feasibility in \eqref{eq:admset} still
gives Theorem~\ref{thm:stab-by-constr}.
The model generates the data values \eqref{eq:targets}.
In these performance statements $J^{\pi_i}$ is the cost
\eqref{eq:Jpi} of holding the $i$th feedback fixed.

\subsection{Polynomial models and sum-of-squares certificates}
\label{subsec:sos}
For multi-indices $\alpha,\beta$ with
nonnegative integer entries, write
$x^\alpha=\prod_{\ell=1}^n x_\ell^{\alpha_\ell}$ and
$|\alpha|=\sum_\ell\alpha_\ell$.
Let each component $g_\ell$ of the dynamics $g$ in~\eqref{eq:sys} be a polynomial as follows.
\begin{equation*}
g_\ell(x,u)=\sum_{|\alpha|+|\beta|\le d_g}
 c_{\ell,\alpha,\beta}x^\alpha u^\beta,\qquad \ell=1,\ldots,n.
\end{equation*}
Let $\phi_d(x)$ collect all monomials $x^\alpha$ with
$|\alpha|\le d$.
The value approximation is parametrized by
\begin{align*}
W_\theta(x)&=\theta^\top\phi_d(x),\qquad N_d=\binom{n+d}{d},\\
\Fset_d&=\{W_\theta:\theta\in\R^{N_d}\}.\nonumber
\end{align*}

Additionally, assume that $V,r$ are polynomial. 
Take $\rho(s)=as$ with $0<a<1$, $b_j(s)=B_js$ with $B_j>0$, and the state and input sets have descriptions $\Xset=\{x:h_a^X(x)\ge0 \forall a\}$, $\Uset=\{u:h_b^U(u)\ge0\ \forall b\}$ using finitely many polynomials.
At a fixed measured $\xi$, Algorithm~\ref{alg:control} solves the polynomial optimization problem
\begin{equation}
\label{eq:polycontrol}
\begin{aligned}
\min_{\bu,x_1,\ldots,x_M}\quad&
 \sum_{j=0}^{M-1}[V(x_j)+r(u_j)]+\gamma W_\theta(x_M)\\
\text{subject to}\quad&
 x_0=\xi,\quad x_{j+1}=g(x_j,u_j),\\
&h_b^U(u_j)\ge0,\quad h_a^X(x_{j+1})\ge0,\\
&B_j^2\xi^\top\xi-x_j^\top x_j\ge0\quad(1\le j<M),\\
&aV(\xi)-V(x_M)\ge0.
\end{aligned}
\end{equation}
Given the target values \eqref{eq:targets}, one possible fitting calculation is $\theta_{i+1}\in\argmin_\theta \sum_{s=1}^{N_s}w_s [\theta^\top\phi_d(\xi^{(s)})-y^{(s)}]^2$, with $w_s>0$.
Eliminate $x_j$ through the polynomial dynamics. Then
$K_c$ is described by finitely many polynomial inequalities
$h_a(\xi,\bu)\ge0$, including $c-V(\xi)\ge0$ and every
constraint of \eqref{eq:polycontrol}. Let
$\Omegac=\{\xi:s_b(\xi)\ge0\ \forall b\}$.
A polynomial $\sigma$ is a sum of squares (SOS) if
$\sigma(z)=\sum_\nu p_\nu(z)^2$. Define $\mathcal Q(h)= \{\sigma_0+\textstyle\sum_a\sigma_a h_a: \sigma_0,\sigma_a\ {\rm SOS}\}$.
Every member of $\mathcal Q(h)$ is nonnegative on the set
defined by $h_a\ge0$.
For fixed polynomial $A,B$ and a polynomial candidate feedback $\widetilde\pi(\xi)=P\phi_{d_\pi}(\xi)\in\R^{mM}$.
For computation, an SOS polynomial has the form $\sigma(z)=\phi_s(z)^\top H\phi_s(z)$ with $H\succeq0$.
Matching coefficients in $p=\sigma_0+\sum_a\sigma_a h_a$ and constraining the matrices $H$ to be positive semidefinite gives a semidefinite feasibility problem.
If we take $A,B,\widetilde\pi$ fixed, these verification problems are convex.

\subsection{Piecewise-affine models and value approximations}
\label{subsec:structured}
Let $\mathcal D_r$ be a finite collection of polyhedra covering
the relevant state--input domain. 
A continuous piecewise-affine
plant has the representation $g(x,u)=A_rx+B_ru+a_r$  for $(x,u)\in\mathcal D_r$.
For a fixed polyhedral partition $\{\mathcal P_l\}_{l=1}^{N_P}$
covering $\Omegac$, use $W_\theta(x)=v_l^\top x+w_l$ for $x\in\mathcal P_l$.
The parameter vector $\theta$ collects all $v_l,w_l$.
Impose equality of the affine formulas on every common face.
This defines the continuous family $\Fset$.
Fitting the values in~\eqref{eq:targets} by least squares with these conditions is a
quadratic optimization problem with linear constraints.
For a simplicial partition, the values at the vertices
equivalently specify the continuous affine interpolants.

At a fixed $\xi$, introduce predicted states and impose, for
every $j$, $\bigvee_r \left[ (x_j,u_j)\in\mathcal D_r,\quad x_{j+1}=A_rx_j+B_ru_j+a_r \right]$.
Binary variables and validated finite bounds encode this finite choice. 
For example, pick
\[
\begin{split}
&\lambda_{jr}\in\{0,1\},\quad \sum_r\lambda_{jr}=1,\\
&H_r[x_j^\top,u_j^\top]^\top
       \le h_r+C_r(1-\lambda_{jr}),\\
&x_{j+1}-A_rx_j-B_ru_j-a_r
       \le D_r(1-\lambda_{jr}),\\
&-x_{j+1}+A_rx_j+B_ru_j+a_r
       \le D_r(1-\lambda_{jr}),
\end{split}
\]
when $\mathcal D_r=\{z:H_rz\le h_r\}$.
The componentwise constants $C_r,D_r$ must dominate possible violations on the bounded variable domain.
The same construction imposes the equality $t=W_\theta(x_M)$ through the regions $\mathcal P_l$.
Moreover, to specify the certificate and cost, one useful case is $V(x)=\max_{a=1,\ldots,N_V}(p_a^\top x+v_a)$, $r(u)=\max_{b=1,\ldots,N_r}(q_b^\top u+r_b)$.

\subsection{Rectified-linear-unit value approximations}
\label{subsec:relu}
A rectified linear unit is $\sigma(s)=\max\{0,s\}$, applied
componentwise to vectors. A value approximation with $L$
hidden layers is
\begin{equation*}
\begin{aligned}
z^0&=x,\\
s^l&=A^lz^{l-1}+a^l,\quad z^l=\sigma(s^l),
                         &&l=1,\ldots,L,\\
W_\theta(x)&=a_{\rm out}^{\top}z^L+b_{\rm out}.
\end{aligned}
\end{equation*}
The parameter vector $\theta$ contains all weights and biases,
and layer sizes are fixed when defining $\Fset$.
This is a continuous piecewise-affine scalar function.
In Algorithm~\ref{alg:valueiteration}, its parameters are fitted
by minimizing
$\sum_s w_s[W_\theta(\xi^{(s)})-y^{(s)}]^2$.
This fitting problem is generally nonconvex. A small achieved
sample loss supplies no uniform certificate.
In Algorithm~\ref{alg:control}, $\theta$ is fixed and the
network is evaluated at the predicted endpoint $x_M$.
For a neuron with preactivation $s$ and validated bounds
$l\le s\le h$, where $l<0<h$, the exact graph $z=\max\{0,s\}$
is $z\ge0$, $z\ge s$, $d\in\{0,1\}$, $z\le hd$, $z\le s-l(1-d)$.
If $h\le0$, set $z=0$. Moreover, if $l\ge0$, set $z=s$.
The two binary cases resulting from $d$ give exactly the inactive and active branches, respectively~\cite{Tjeng.2019}.

\section{Robustification again Model Uncertainty}
\label{sec:robust}
Let  $\widehat g$ denote an approximation of the true dynamics $g$.
Given an fs-CLF $V$ of length $M$, $\mathcal{K}$-bounds $b_i$ as in~\eqref{eq:K_bounded_g_closed_loop}, and the approximate dynamics $\widehat g$, we aim to design a learned finite-step control law.
Towards this end, we show the closeness of the resulting system trajectories from the true system dynamics $g$ and the learned one $\widehat g$. This is obtained by making an assumption on the mismatch between the approximate model and the true one as follows.

Given a control sequence $\bu$ and $\widehat x_0=\xi$, the learned system is represented by
\begin{equation}
\label{eq:modelprediction}
\widehat x_{j+1}=\widehat g(\widehat x_j,u_j).
\end{equation}

\begin{assumption}
\label{ass:modelinformation}
Pick a compact set $\mathcal D\subseteq\Xset$ containing $\Omegac$ as in~\eqref{eq:omegac}.
Assume that the learned dynamics $\widehat g:\Rn\times\Rm\to\Rn$ as in~\eqref{eq:modelprediction} is continuous.
There exists a nonnegative function $d : \mathcal D \times \Uset \to \Rpos$ such that for all $(x,u)\in\mathcal D\times\Uset$ the following holds
\begin{equation}
\label{eq:modelerror}
\norm{g(x,u)-\widehat g(x,u)}\le d(x,u) .
\end{equation}
Moreover, there exists a known $L_V\in\Rpos$ satisfying
\begin{equation}\label{eq:V_incremental_bound}
|V(x)-V(y)|\le L_V \norm{x-y} \qquad \forall x,y\in\mathcal D.
\end{equation}
\end{assumption}

\begin{assumption}\label{ass:standing_learned_syst}
Let $\mathcal{D} \subseteq \Xset$ as in Assumption~\ref{ass:modelinformation} be given.
Suppose that the learned dynamics~\eqref{eq:modelprediction} satisfies Assumption~\ref{ass:standing} for some $\Uset$ and with $\mathcal{D}$ in place of $\Xset$.
With slight abuse of notation, let $V$ be a fs-CLF for~\eqref{eq:modelprediction}. 
\end{assumption}

In line of Assumption~\ref{ass:standing_learned_syst}, denote the associated finite-step admissible feedback set by 
\begin{equation}
\label{eq:practicaladmissible}
\begin{gathered}
\AdmL(\xi):=\Big\{\bu\in\Uset^M: \widehat{g}^j(\xi,\bu)\in\mathcal{D} \subseteq \Xset\ (j=1,\ldots,M),\\
 \norm{\widehat x_j}\le b_j(\norm{\xi})\ (1\le j<M),\\
V(\widehat x_M)\le\rho(V(\xi))\Big\}.
\end{gathered}
\end{equation}
Let $L_g\ge0$ be a known constant satisfying
\begin{equation}
\label{eq:truestateL}
\norm{g(\xi_1,\mu)-g(\xi_2,\mu)}\le L_g\norm{\xi_1-\xi_2},
\,\, \xi_1,\xi_2\in\mathcal D,\, \mu\in\Uset.
\end{equation}
We provide forward invariance and practical asymptotic stability certificate as follows.

\begin{theorem}
\label{thm:practical}
Let Assumptions~\ref{ass:modelinformation}
and~\ref{ass:standing_learned_syst} hold with
$\Omegac\subseteq\mathcal D\subseteq\Xset$ for some $c>0$.
Let $\alpha_1,\alpha_2$ satisfy~\eqref{eq:properV} on $\mathcal D$.
Suppose the true dynamics $g$ satisfies~\eqref{eq:truestateL}.
Assume that there exists a known $\bar d\ge0$ satisfying
\begin{equation}
\label{eq:uniformmodelerror}
d(x,u)\le\bar d \qquad\forall(x,u)\in\mathcal D\times\Uset.
\end{equation}
Define
\begin{equation}
\label{eq:uniformpredictionerrors}
\bar\Delta_j:=\bar d\sum_{l=0}^{j-1}L_g^l 0\le j\le M,\qquad \Delta_0:=L_V\bar\Delta_M,
\end{equation}
where an empty sum is zero.
For $\xi\in\Omegac$ and $\bu\in\AdmL(\xi)$, define the
prediction-error bounds and the associated safe subset by
\begin{subequations}\label{eq:nominaldomainmargin}
\begin{equation}\label{eq:sequencepredictionerrors}
\begin{aligned}
e_0(\xi,\bu)&:=0,\\
e_{j+1}(\xi,\bu)&:=L_g e_j(\xi,\bu)+d(\widehat x_j,u_j), \quad 0\le j<M.
\end{aligned}
\end{equation}
\begin{equation}\label{eq:safeadmissible}
\begin{aligned}
\mathcal S_M(\xi):=\big\{\bu\in\AdmL(\xi):{}&
\Bball(\widehat x_j,e_j(\xi,\bu))\subseteq\mathcal D\\
&\text{for }j=1,\ldots,M\big\}.
\end{aligned}
\end{equation}
Suppose
\begin{equation}\label{eq:safenonempty}
\mathcal S_M(\xi)\ne\varnothing
\qquad\forall\xi\in\Omegac.
\end{equation}
\end{subequations}
By~\eqref{eq:uniformmodelerror} and~\eqref{eq:sequencepredictionerrors},
\[
e_j(\xi,\bu)
=\sum_{l=0}^{j-1}L_g^{j-1-l}d(\widehat x_l,u_l)
\le\bar\Delta_j,\qquad 0\le j\le M.
\]
Assume
\begin{equation}
\label{eq:uniformmismatchsmallness}
\Delta_0<c-\rho(c).
\end{equation}
Set
\begin{equation}
\label{eq:practicalratechoice}
\begin{gathered}
\vartheta:=\frac12\left(1+\frac{\Delta_0}{c-\rho(c)}\right),\\
\tilde\rho(s):=(1-\vartheta)\rho(s)+\vartheta s.
\end{gathered}
\end{equation}
\begin{equation}
\label{eq:practicalthreshold}
a:=\max\left\{s\in[0,c]:s-\rho(s)\le
                    \frac{\Delta_0}{\vartheta}\right\}.
\end{equation}
Then $\tilde\rho\in\Kinf$, $\tilde\rho<\id$, and $0\le a<c$.

For $\bx(0)=\xi\in\Omegac$, choose at each $i\in\N$
\begin{equation}
\label{eq:practicalexecutedblock}
\bu_i\in\mathcal S_M(\bx(i)),\qquad
\bx(i+1)=\Fm(\bx(i),\bu_i).
\end{equation}
Then the resulting true trajectory satisfies
$x(k)\in\mathcal D$ for all $k\in\N$, and
$\Omegac$ is invariant at block boundaries.
Moreover,
\begin{equation}
\label{eq:practicalcomparison}
V(\bx(i))\le\max\{\tilde\rho^i(V(\xi)),a\}
\qquad\forall i\in\N,
\end{equation}
where $\tilde\rho^i$ denotes $i$-fold composition.
Set $b_0=\id$.
For $i\in\N$ and $0\le j<M$,
\begin{equation}
\label{eq:intermiateTrajEstimate}
\begin{aligned}
\norm{x(iM+j)}\le\max\big\{&
 b_j\circ\alpha_1^{-1}\circ\tilde\rho^i(V(\xi)),\\
&b_j(\alpha_1^{-1}(a))\big\}+\bar\Delta_j.
\end{aligned}
\end{equation}
Define
\begin{equation}
\label{eq:practicalradius}
D:=\max_{0\le j<M}
       \{b_j(\alpha_1^{-1}(a))+\bar\Delta_j\}.
\end{equation}
There exists $\beta\in\KL$ such that
\begin{equation}
\label{eq:practicalAsympticEstimate}
\norm{x(k)}\le\beta(\norm{\xi},k)+D
\qquad\forall k\in\N,\quad\xi\in\Omegac.
\end{equation}
\end{theorem}

\begin{proof}
Pick $\xi\in\Omegac$ and
$\bu_0=(u_0,\ldots,u_{M-1})\in\mathcal S_M(\xi)\subseteq\AdmL(\xi)$.
Write $e_j=e_j(\xi,\bu_0)$ as in~\eqref{eq:sequencepredictionerrors}.
Set $x_j=g^j(\xi,\bu_0)$ and
$\widehat x_j=\widehat g^j(\xi,\bu_0)$ for $0\le j\le M$.
By~\eqref{eq:practicaladmissible}, $\widehat x_j\in\mathcal D$.
If $x_j\in\mathcal D$, equations~\eqref{eq:modelerror}
and~\eqref{eq:truestateL} give
\begin{align}
\label{eq:truepredictionrecursion}
\norm{x_{j+1}-\widehat x_{j+1}}
&=\norm{g(x_j,u_j)-\widehat g(\widehat x_j,u_j)}\notag\\
&\le\norm{g(x_j,u_j)-g(\widehat x_j,u_j)}\notag\\
&\quad+\norm{g(\widehat x_j,u_j)-\widehat g(\widehat x_j,u_j)}\notag\\
&\le L_g\norm{x_j-\widehat x_j}+d(\widehat x_j,u_j).
\end{align}
We prove by induction on $j$ that
\begin{equation}
\label{eq:truestatedomainerror}
\begin{gathered}
x_j\in\mathcal D,\qquad
\norm{x_j-\widehat x_j}\le e_j,\\
e_j=\sum_{l=0}^{j-1}L_g^{j-1-l}d(\widehat x_l,u_l)
\le\bar\Delta_j,\qquad 0\le j\le M.
\end{gathered}
\end{equation}
For $j=0$, $x_0=\widehat x_0=\xi\in\mathcal D$ and $e_0=0$.
Suppose~\eqref{eq:truestatedomainerror} holds for some $0\le j<M$.
Since $x_j,\widehat x_j\in\mathcal D$, equations
\eqref{eq:truepredictionrecursion} and~\eqref{eq:sequencepredictionerrors}
give
\begin{align*}
\norm{x_{j+1}-\widehat x_{j+1}}
&\le L_g e_j+d(\widehat x_j,u_j)\\
&=e_{j+1}
 =\sum_{l=0}^{j}L_g^{j-l}d(\widehat x_l,u_l)\\
&\le\bar d\sum_{l=0}^{j}L_g^{j-l}
 =\bar\Delta_{j+1}.
\end{align*}
Since $\norm{x_{j+1}-\widehat x_{j+1}}\le e_{j+1}$ from the above inequality, we have $x_{j+1}\in\Bball(\widehat x_{j+1},e_{j+1})$.
Moreover, $\bu_0\in\mathcal S_M(\xi)$ implies
$\Bball(\widehat x_{j+1},e_{j+1})\subseteq\mathcal D$
by~\eqref{eq:safeadmissible}.
Hence, $x_{j+1}\in\mathcal D$.
This proves~\eqref{eq:truestatedomainerror} by induction, including
the case $e_{j+1}=0$.
In particular, $x_M,\widehat x_M\in\mathcal D$.
From~\eqref{eq:V_incremental_bound},
\eqref{eq:practicaladmissible}, and~\eqref{eq:truestatedomainerror},
\begin{align}
\label{eq:trueendpointperturbation}
V(\Fm(\xi,\bu_0))
&\le V(\widehat x_M)+|V(x_M)-V(\widehat x_M)|\notag\\
&\le\rho(V(\xi))+L_V\norm{x_M-\widehat x_M}\notag\\
&\le\rho(V(\xi))+
 L_V\sum_{l=0}^{M-1}L_g^{M-1-l}d(\widehat x_l,u_l)\notag\\
&\le\rho(V(\xi))+L_V\bar d\sum_{l=0}^{M-1}L_g^l\notag\\
&=\rho(V(\xi))+\Delta_0.
\end{align}
Since $V(\xi)\le c$, equations~\eqref{eq:uniformmismatchsmallness}
and~\eqref{eq:trueendpointperturbation} yield $V(x_M)\le\rho(c)+\Delta_0<c$.
Hence $x(M)\in\Omegac$.
For any $i\in\N$ with $x(iM)\in\Omegac$, choose
$\bu_i\in\mathcal S_M(x(iM))$, which exists
by~\eqref{eq:safenonempty}, and repeat
\eqref{eq:truepredictionrecursion}--\eqref{eq:trueendpointperturbation}
with initial state $x(iM)$.
It follows that
\begin{equation}
\label{eq:allblocksdomain}
\begin{gathered}
x(iM+j)\in\mathcal D\quad(0\le j\le M),\\
V(\bx(i+1))\le\rho(V(\bx(i)))+\Delta_0<c.
\end{gathered}
\end{equation}
Induction on $i$ proves $\bx(i)\in\Omegac$ for every $i\in\N$
and $x(k)\in\mathcal D$ for every $k\in\N$.

To prove~\eqref{eq:practicalcomparison}, we use the threshold construction in~\cite[Theorem~4.1]{Geiselhart.2016} and~\cite[Theorem~7]{Noroozi.2018} as follows.
Equations~\eqref{eq:uniformmismatchsmallness} and~\eqref{eq:practicalratechoice} imply $0\le\frac{\Delta_0}{c-\rho(c)}<\vartheta<1$, $\frac{\Delta_0}{\vartheta}<c-\rho(c)$.
The set in~\eqref{eq:practicalthreshold} is nonempty because it contains zero.
It is closed in $[0,c]$ by continuity of $\rho$, and it does not contain $c$.
Therefore its maximum $a$ exists and $0\le a<c$.
For $a<s\le c$, maximality gives $s-\rho(s)>\Delta_0/\vartheta$.
Letting $s\downarrow a$ and using the defining inequality at $a$ gives
\begin{equation}
\label{eq:thresholdgap}
a-\rho(a)=\frac{\Delta_0}{\vartheta},\qquad
s-\rho(s)\ge\frac{\Delta_0}{\vartheta}\quad(a\le s\le c).
\end{equation}
The function $\tilde\rho$ in~\eqref{eq:practicalratechoice} is continuous,
$\tilde\rho(0)=0$, and, for $0\le s<t$, $\tilde\rho(t)-\tilde\rho(s)=(1-\vartheta)(\rho(t)-\rho(s))+\vartheta(t-s)>0$, $\tilde\rho(s)\ge\vartheta s$, $\tilde\rho(s) < (1-\vartheta)s+\vartheta s=s\quad(s>0)$.
Hence $\tilde\rho\in\Kinf$ and $\tilde\rho<\id$.

For $a\le s\le c$, equation~\eqref{eq:thresholdgap} implies $\rho(s)+\Delta_0
\le\rho(s)+\vartheta(s-\rho(s)) =\tilde\rho(s)$.
For $0\le s\le a$, monotonicity of $\rho$ gives
\begin{align*}
\rho(s)+\Delta_0
&\le\rho(a)+\Delta_0\\
&=a-\frac{\Delta_0}{\vartheta}+\Delta_0
 =a-\frac{1-\vartheta}{\vartheta}\Delta_0\le a.
\end{align*}
Combining these two inequalities with~\eqref{eq:allblocksdomain} yields
\begin{equation}
\label{eq:MstepDecay}
V(\bx(i+1))\le\max\{\tilde\rho(V(\bx(i))),a\}
\qquad\forall i\in\N.
\end{equation}
At $i=0$, inequality~\eqref{eq:practicalcomparison} follows from
$V(\bx(0))=V(\xi)$.
If it holds at $i$, then~\eqref{eq:MstepDecay} and monotonicity of
$\tilde\rho$ give
\begin{align*}
V(\bx(i+1))
&\le\max\{\tilde\rho(\max\{\tilde\rho^i(V(\xi)),a\}),a\}\\
&=\max\{\tilde\rho^{i+1}(V(\xi)),\tilde\rho(a),a\}\\
&=\max\{\tilde\rho^{i+1}(V(\xi)),a\}.
\end{align*}
This proves~\eqref{eq:practicalcomparison} by induction.

Now take $k=iM+j$, where $i\in\N$ and $0\le j<M$.
Use $b_0=\id$ and $\bar\Delta_0=0$.
Equation~\eqref{eq:truestatedomainerror}, applied to the block starting
at $x(iM)$, gives
\begin{align*}
\norm{x(iM&+j)}
\le\norm{g^j(\bx(i),\bu_i)-\widehat g^j(\bx(i),\bu_i)}+\norm{\widehat g^j(\bx(i),\bu_i)}\\
&\le\bar\Delta_j+b_j(\norm{\bx(i)})\\
&\le\bar\Delta_j+b_j\circ\alpha_1^{-1}(V(\bx(i)))\\
&\le\bar\Delta_j+b_j\circ\alpha_1^{-1}
             (\max\{\tilde\rho^i(V(\xi)),a\})\\
&=\bar\Delta_j+\max\left\{
b_j\circ\alpha_1^{-1}\circ\tilde\rho^i(V(\xi)), b_j(\alpha_1^{-1}(a))\right\}.
\end{align*}
This proves~\eqref{eq:intermiateTrajEstimate}.

For $r\ge0$ and $i\in\N$, define
\begin{equation}
\label{eq:betablock}
\beta_i(r):=\max_{0\le j<M}
      b_j\circ\alpha_1^{-1}\circ\tilde\rho^i\circ\alpha_2(r).
\end{equation}
Then $\beta_i\in\Kinf$ and $\beta_{i+1}(r)\le\beta_i(r)$.
For each $s\ge0$, the sequence $\tilde\rho^i(s)$ is nonincreasing and
bounded below by zero.
Its limit satisfies $s_\infty:=\lim_{i\to\infty}\tilde\rho^i(s) =\tilde\rho(s_\infty), \tilde\rho<\id$ which implies $s_\infty=0$.
Consequently, $\beta_i(r)\to0$ for each $r\ge0$.
Define a continuous interpolation by
\begin{equation}
\label{eq:betainterpolation}
\begin{gathered}
\beta(r,t)=\beta_0(r)\qquad(0\le t\le M),\\
\beta(r,(i+\theta)M)
 =(1-\theta)\beta_{i-1}(r)+\theta\beta_i(r),\\
i\in\N^*,\qquad0\le\theta\le1.
\end{gathered}
\end{equation}
The two formulas agree at $t=M$, and consecutive interpolation
intervals agree at their common endpoints.
Since $\beta_i\in\Kinf$ and $\beta_{i+1}\le\beta_i$, one can easily see that $\beta\in\KL$.

For $k=iM+j$, $0\le j<M$, equation~\eqref{eq:betainterpolation} gives
$\beta(\norm{\xi},k)\ge\beta_i(\norm{\xi})$.
Finally, equations~\eqref{eq:properV},
\eqref{eq:intermiateTrajEstimate}, and~\eqref{eq:practicalradius} imply $\norm{x(iM+j)}\le\beta_i(\norm{\xi})+b_j(\alpha_1^{-1}(a))+\bar\Delta_j \le\beta(\norm{\xi},iM+j)+D$,
which proves~\eqref{eq:practicalAsympticEstimate}.
\end{proof}

For a fixed $\widehat g$ as in~\eqref{eq:modelprediction}, define the predicted block cost
and the numerical control objective by
\begin{align}
\label{eq:modelblockcost}
\widehat L(\xi,\bu)&:=\sum_{j=0}^{M-1}
                       [V(\widehat x_j)+r(u_j)],\\
\label{eq:modelcontrolobjective}
\widehat Q_{\Vhat}(\xi,\bu)&:=\widehat L(\xi,\bu)
                                      +\gamma\Vhat(\widehat x_M).
\end{align}
Now we provide the counterpart of Algorithm~\ref{alg:control}, where the learned dynamics $\widehat g$ is used in place of the true dynamics $g$.

\begin{controlalgorithm}{Feedback control with certified model uncertainty}
\label{alg:robustcontrol}
\textit{Input Data:} $\Xset,\Uset,V,\rho,M,b_j$ from
\eqref{eq:admset}, $\Omegac\subseteq\mathcal D$ from
\eqref{eq:omegac}, $\xi=x(kM)\in\Omegac$,
$\ell,r,\gamma$ from \eqref{eq:stagecost}--\eqref{eq:Jpi},
$\Vhat\in C(\Omegac)$, a learned model $\widehat g$
satisfying Assumption~\ref{ass:standing_learned_syst},
$d,L_V,\mathcal D$ from Assumption~\ref{ass:modelinformation},
$L_g$ from~\eqref{eq:truestateL} and $\bar d\ge0$
satisfying~\eqref{eq:uniformmodelerror} and~\eqref{eq:uniformmismatchsmallness},
the sets~\eqref{eq:safeadmissible} satisfying~\eqref{eq:safenonempty},
and a numerical solver for the constrained problems in steps 2--3.
\begin{enumerate}[label=\arabic*.,leftmargin=1.8em]
\item Introduce the optimization variables
$\bu=(u_0,\ldots,u_{M-1})\in\Uset^M$.
Set $\widehat x_0=\xi$ and form the predictions
\eqref{eq:modelprediction}.
Form the bounds in~\eqref{eq:sequencepredictionerrors}.
\item Compute a feasible solution $\bu^0$ of
\begin{equation}
\label{eq:learnedfeasibility}
\min_{\bu\in\Uset^M}0
\quad\text{subject to}\quad
\bu\in\mathcal S_M(\xi)
\end{equation}
using~\eqref{eq:practicaladmissible} and~\eqref{eq:safeadmissible}.
If no feasible solution is obtained, return failure.
\item Minimize $\widehat Q_{\Vhat}(\xi,\bu)$ in
\eqref{eq:modelcontrolobjective} over $\bu\in\mathcal S_M(\xi)$, starting from $\bu^0$.
For a computed solution $\bu^1$, set
$\bu^{\mathrm f}=\bu^1$ if $\bu^1\in\mathcal S_M(\xi)$ and $\widehat Q_{\Vhat}(\xi,\bu^1)\le
 \widehat Q_{\Vhat}(\xi,\bu^0)$.
Otherwise set $\bu^{\mathrm f}=\bu^0$.
\item Return $\bu^{\mathrm f}\in\mathcal S_M(\xi)$ and $(\widehat x_j(\xi,\bu^{\mathrm f}))_{j=0}^M$.
\end{enumerate}
\end{controlalgorithm}

\section{Model Identification and Constraint-Preserving Data Collection}
\label{sec:identification}
Algorithm~\ref{alg:robustcontrol} relies on the availability of learned dynamics $\widehat g$ satisfying Assumption~\ref{ass:modelinformation}.
This section constructs such a function $d$ required by Assumption~\ref{ass:modelinformation} from measured states.
The construction assumes known bounds on measurement errors and on the variation of $g$, while its inequalities hold for every error consistent with those bounds, without a probability qualification.

\subsection{A deterministic error bound from input-state measurements}
Let $z=(x^\top,u^\top)^\top\in\mathcal D\times\Uset$.
For $t\ge1$, let the joint input-state measured data be given by
\begin{equation}
\label{eq:transitiondata}
\begin{gathered}
\mathcal D_t^{\rm obs}=\{(z_s,y_s,\nu_s):s=1,\ldots,t\},\\
z_s\in\mathcal D\times\Uset,\quad
\norm{y_s-g(z_s)}\le\nu_s,\quad\nu_s\ge0 ,
\end{gathered}
\end{equation}
where $g(z_s)$ denotes $g(x_s,u_s)$, the location $z_s$ is known, $y_s$ is the measured next state, and $\nu_s$
is a known measurement-error bound.
Assume that a joint state--input Lipschitz factor is given by
\begin{equation}
\label{eq:jointplantL}
\norm{g(z)-g(z')}\le L_{\rm sm}\norm{z-z'},
\qquad z,z'\in\mathcal D\times\Uset .
\end{equation}
Let $\widehat g$ be a fitted approximation of $g$.
If the model is fitted to the measured transition data $\mathcal D_t^{\rm obs}$, it may depend on the data index $t$.
For simplicity, we omit this dependence in the notation for $\widehat g$. When comparing the model-error bounds at different data indices in the following proposition, the model $\widehat g$ is kept fixed.

\begin{proposition}
\label{prop:setmembership}
Let the measured transition data satisfy
\eqref{eq:transitiondata} and assume
\eqref{eq:jointplantL}.
For $t\ge1$, the function
\begin{equation}
\label{eq:smerror}
d_t(z):=\min_{1\le s\le t}
\big\{\norm{\widehat g(z)-y_s}
+\nu_s+L_{\rm sm}\norm{z-z_s}\big\}
\end{equation}
satisfies~\eqref{eq:modelerror}.
For a fixed model $\widehat g$, if the additional measured transition satisfies \eqref{eq:transitiondata}, then $d_{t+1}(z)\le d_t(z)$ for all $z\in\mathcal D\times\Uset$.
Moreover, if $\nu_s=0$ for all $s\ge1$ and $\lim_{t\to\infty}\min_{1\le s\le t}\norm{z-z_s}=0$
for all $z\in\mathcal D\times\Uset$, then $\lim_{t\to\infty}d_t(z) =\norm{\widehat g(z)-g(z)}$
for all $z\in\mathcal D\times\Uset$.
\end{proposition}
\begin{proof}
Fix $z\in\mathcal D\times\Uset$.
For each measured transition indexed by $s\in\{1,\ldots,t\}$, equations \eqref{eq:transitiondata} and
\eqref{eq:jointplantL} yield
\begin{align*}
\norm{\widehat g(z)-g(z)}
&\le \norm{\widehat g(z)-y_s}
    +\norm{y_s-g(z_s)}\\
&\quad+\norm{g(z_s)-g(z)}\\
&\le \norm{\widehat g(z)-y_s}
    +\nu_s+L_{\rm sm}\norm{z-z_s}.
\end{align*}
Taking the minimum over $s\in\{1,\ldots,t\}$ gives $\norm{\widehat g(z)-g(z)}\le d_t(z)$.
For a fixed model $\widehat g$, it follows from~\eqref{eq:smerror} that
\begin{align*}
d_{t+1}(z)=\min\big\{&d_t(z), \norm{\widehat g(z)-y_{t+1}}\\
&+\nu_{t+1}+L_{\rm sm}\norm{z-z_{t+1}}\big\} \le d_t(z).
\end{align*}

Now suppose that $\nu_s=0$ for all $s\ge1$.
By \eqref{eq:transitiondata}, we have $y_s=g(z_s)$.
It follows from \eqref{eq:jointplantL} that
\begin{align*}
\norm{\widehat g(z)-y_s}
&\le \norm{\widehat g(z)-g(z)}
    +\norm{g(z)-g(z_s)}\\
&\le \norm{\widehat g(z)-g(z)}
    +L_{\rm sm}\norm{z-z_s}.
\end{align*}
Substituting this inequality into \eqref{eq:smerror}
yields $d_t(z)\le\norm{\widehat g(z)-g(z)} +2L_{\rm sm}\min_{1\le s\le t}\norm{z-z_s}$.
Combining this inequality with the model-error bound proved above gives $0\le d_t(z)-\norm{\widehat g(z)-g(z)} \le 2L_{\rm sm}\min_{1\le s\le t}\norm{z-z_s}$.
By the assumed limit condition, the right-hand side converges to zero as $t\to\infty$. This proves $\lim_{t\to\infty}d_t(z) =\norm{\widehat g(z)-g(z)}$.
Since $z\in\mathcal D\times\Uset$ was arbitrary, the proof is complete.
\end{proof}

If the model $\widehat g$ is refitted as additional measured transitions become available, the term
$\norm{\widehat g(z)-y_s}$ in \eqref{eq:smerror} may change, and the inequality $d_{t+1}(z)\le d_t(z)$
is no longer guaranteed.

To retain the model-error information independently of the choice of the fitted model, construct
\emph{uncertainty} sets $\mathcal H_t(z)$ containing the possible values of $g(z)$ for every
$z\in\mathcal D\times\Uset$.
For the measured data \eqref{eq:transitiondata}, define
\begin{equation}
\label{eq:nesteduncertainty}
\mathcal H_t(z)=\mathcal H_0(z)\cap \bigcap_{s=1}^t \Bball(y_s,\nu_s+L_{\rm sm}\norm{z-z_s}),
\end{equation}
where $\mathcal H_0(z)$ is a known compact set containing
$g(z)$. Equations \eqref{eq:transitiondata} and
\eqref{eq:jointplantL} imply that each intersected ball
contains $g(z)$.
By construction, the set $\mathcal H_t$ is nonempty and compact, and
\begin{equation}
\label{eq:nestedinformationproperty}
g(z)\in\mathcal H_{t+1}(z)\subseteq\mathcal H_t(z)
\qquad\forall z\in\mathcal D\times\Uset.
\end{equation}
For any continuous model $h:\Rn\times\Rm\to\Rn$, define
\begin{equation}
\label{eq:retainedmodelerror}
d^{\mathcal H_t}(z,h):= \max_{y\in\mathcal H_t(z)}\norm{y-h(z)}.
\end{equation}
Compactness gives a finite maximum at each $z$.
Since $g(z)\in\mathcal H_t(z)$, we have
\begin{equation}
\label{eq:retainederrorvalidity}
\norm{g(z)-h(z)}\le d^{\mathcal H_t}(z,h)
\qquad\forall z\in\mathcal D\times\Uset.
\end{equation}
Thus \eqref{eq:retainedmodelerror} gives the bound required
by \eqref{eq:modelerror} for any chosen model $h$.
As in \eqref{eq:modelprediction}, denote $\widehat x_0^h=\xi$, $\widehat x_{j+1}^h=h(\widehat x_j^h,u_j)$, $0\le j<M$.
Define the learned-model admissible set by
\begin{subequations}\label{eq:retainedadmissibility}
\begin{equation}
\label{eq:setadmissibility}
\begin{gathered}
\mathcal U_{M}^{h}(\xi):=
\Big\{\bu\in\Uset^M:\widehat x_j^h\in\mathcal D\ (1\le j\le M),\\
\norm{\widehat x_j^h}\le b_j(\norm{\xi})\ (1\le j<M),\\
V(\widehat x_M^h)\le\rho(V(\xi))\Big\}.
\end{gathered}
\end{equation}
For $\bu\in\mathcal U_{M}^{h}(\xi)$, define
\begin{equation}\label{eq:retainedpredictionerrors}
\begin{aligned}
e_{0,t}^h(\xi,\bu)&:=0,\\
e_{j+1,t}^h(\xi,\bu)&:=L_g e_{j,t}^h(\xi,\bu)
 +d^{\mathcal H_t}((\widehat x_j^h,u_j),h),\\
&\hspace{5em}0\le j<M,
\end{aligned}
\end{equation}
and
\begin{align}
\mathcal S_{M,t}^{\mathcal H}&(\xi,h)  := \nonumber\\
& \big\{\bu\in\mathcal U_{M}^{h}(\xi): \Bball(\widehat x_j^h,e_{j,t}^h(\xi,\bu))
 \subseteq\mathcal D, \,  1\le j\le M \big\}. \label{eq:retainedsafeadmissibility}
\end{align}

\end{subequations}
For a fixed learned model \(h\), the nominal admissible set \(\mathcal U_M^h(\xi)\) is independent of the measurement index \(t\) as its defining constraints depend only on the predicted trajectories generated by \(h\). The uncertainty set \(\mathcal H_t\), however, determines the prediction-error bounds in \(\eqref{eq:retainedpredictionerrors}\) and therefore affects the robustly safe admissible set in \(\eqref{eq:retainedsafeadmissibility}\). As shown next, additional valid measurements can reduce these bounds and enlarge the robustly safe admissible set while \(h\) remains fixed.

The following result ensures preservation of feasible safe
sequences for a fixed learned model $h$.
\begin{proposition}
\label{prop:nestedfeasibility}
Fix $\mathcal D,V,M,\rho,b_j$, $L_g\ge0$, and a model $h$.
If \eqref{eq:nestedinformationproperty} holds, then
\begin{equation}
\label{eq:retainederrormonotonicity}
d^{\mathcal H_{t+1}}(z,h)\le d^{\mathcal H_t}(z,h)
\qquad\forall z\in\mathcal D\times\Uset,
\end{equation}
and
\begin{equation}\label{eq:nestedsafefeasible}
\mathcal S_{M,t}^{\mathcal H}(\xi,h) \subseteq\mathcal S_{M,t+1}^{\mathcal H}(\xi,h)
\qquad\forall\xi\in\Omegac.
\end{equation}
\end{proposition}
\begin{proof}
For every $z\in\mathcal D\times\Uset$, \eqref{eq:nestedinformationproperty} gives
\begin{align*}
d_{t+1}^{\mathcal H}(z,h)
&=\max_{y\in\mathcal H_{t+1}(z)}\norm{y-h(z)}\\
&\le\max_{y\in\mathcal H_t(z)}\norm{y-h(z)}
 =d^{\mathcal H_t}(z,h),
\end{align*}
which proves~\eqref{eq:retainederrormonotonicity}.
Fix $\xi,\bu,h$.
For $\bu\in\mathcal S_{M,t}^{\mathcal H}(\xi,h)$,
the nominal predictions at indices $t$ and $t+1$ coincide.
Equations~\eqref{eq:retainederrormonotonicity}
and~\eqref{eq:retainedpredictionerrors}, with $L_g\ge0$,
give by induction $e_{j,t+1}^h(\xi,\bu)\le e_{j,t}^h(\xi,\bu)$, for $0\le j\le M$.
Hence, for every $1\le j\le M$,
$\Bball(\widehat x_j^h,e_{j,t+1}^h(\xi,\bu)) \subseteq\Bball(\widehat x_j^h,e_{j,t}^h(\xi,\bu))
\subseteq\mathcal D$.
This proves $\bu\in\mathcal S_{M,t+1}^{\mathcal H}(\xi,h)$ and therefore~\eqref{eq:nestedsafefeasible}.
\end{proof}

\subsection{Closed-loop identification with retained certificates}

Proposition~\ref{prop:nestedfeasibility} preserves the nominal set
\eqref{eq:setadmissibility}, enlarges the safe subset
\eqref{eq:retainedsafeadmissibility}, and decreases the error bound
\eqref{eq:retainedmodelerror} when $h$ is fixed.

For a refitted model, the following conditions preserve
the estimates of Theorem~\ref{thm:practical}.

\begin{proposition}
\label{prop:balltoset}
Fix $\mathcal D,\Uset,V,M,\rho,b_j,c,L_g,L_V$, with
$\Omegac\subseteq\mathcal D\subseteq\Xset$ as in~\eqref{eq:omegac}.
Let $\alpha_1,\alpha_2$ satisfy~\eqref{eq:properV} on $\mathcal D$,
and assume~\eqref{eq:V_incremental_bound},
\eqref{eq:truestateL}, and~\eqref{eq:nestedinformationproperty}.
Fix $\bar d\ge0$ satisfying~\eqref{eq:uniformmismatchsmallness}, with
$\bar\Delta_j,\Delta_0$ defined by~\eqref{eq:uniformpredictionerrors}.
At an initial data index $t_0$, let a model $h_0$
satisfy Assumption~\ref{ass:standing_learned_syst}
with the fixed $V,M,\rho,b_j$ and the conditions
\begin{equation}
\label{eq:modelupdateconditions}
\begin{gathered}
\mathcal S_{M,t}^{\mathcal H}(\xi,h)\ne\varnothing
                   \qquad\forall\xi\in\Omegac,\\
d^{\mathcal H_t}(z,h)\le\bar d
                   \qquad\forall z\in\mathcal D\times\Uset.
\end{gathered}
\end{equation}
at $(t,h)=(t_0,h_0)$.
For each subsequent block index $i\ge1$, set $t_i=t_0+iM$.
Use a new model $h_i$ only if it satisfies Assumption~\ref{ass:standing_learned_syst} with the same
$V,M,\rho,b_j$ and~\eqref{eq:modelupdateconditions} at $(t,h)=(t_i,h_i)$.
Otherwise set $h_i=h_{i-1}$.
Then~\eqref{eq:modelupdateconditions} holds at every $(t_i,h_i)$.
In particular, the following holds
\begin{equation}
\label{eq:identifiedcontrol}
\mathcal S_{M,t_i}^{\mathcal H}(\xi,h_i)\ne\varnothing\qquad\forall\xi\in\Omegac,\quad i\in\N,
\end{equation}
and
\begin{equation}
\label{eq:identifieduniformerror}
\norm{g(z)-h_i(z)}\le d_{t_i}^{\mathcal H}(z,h_i)\le\bar d
\qquad z\in\mathcal D\times\Uset,\quad i\in\N.
\end{equation}
For $x(0)=\xi\in\Omegac$, apply any sequence satisfying
\begin{equation}
\label{eq:identifiedexecution}
\bu_i\in\mathcal S_{M,t_i}^{\mathcal H}(\bx(i),h_i),\qquad
\bx(i+1)=\Fm(\bx(i),\bu_i).
\end{equation}
Then $x(k)\in\mathcal D$ for every $k\in\N$,
$\bx(i)\in\Omegac$ for every $i\in\N$, and
\eqref{eq:practicalcomparison}, \eqref{eq:intermiateTrajEstimate},
and~\eqref{eq:practicalAsympticEstimate} hold with the same
$\tilde\rho,a,D,\beta$ as in Theorem~\ref{thm:practical}.
\end{proposition}
\begin{proof}
The conditions~\eqref{eq:modelupdateconditions} hold at $(t_0,h_0)$
by assumption.
Suppose they hold at $(t_{i-1},h_{i-1})$.
If $h_i=h_{i-1}$, equations~\eqref{eq:retainederrormonotonicity}
and~\eqref{eq:nestedsafefeasible}, applied successively from $t_{i-1}$ to $t_i$,
give $d_{t_i}^{\mathcal H}(z,h_i) \le d_{t_{i-1}}^{\mathcal H}(z,h_{i-1})\le\bar d$ for all $z\in\mathcal D\times\Uset$, and $\varnothing\ne\mathcal S_{M,t_{i-1}}^{\mathcal H}(\xi,h_{i-1}) \subseteq\mathcal S_{M,t_i}^{\mathcal H}(\xi,h_i)$ for all $\xi\in\Omegac$.
Thus both conditions in~\eqref{eq:modelupdateconditions}
are preserved when the previous model is retained.
If $h_i\ne h_{i-1}$, all inequalities in~\eqref{eq:modelupdateconditions}
hold by the model-selection rule.
Induction proves~\eqref{eq:modelupdateconditions} at every $(t_i,h_i)$,
including~\eqref{eq:identifiedcontrol}.
Equation~\eqref{eq:retainederrorvalidity} gives~\eqref{eq:identifieduniformerror}.

Pick $\bx(i)\in\Omegac$ and a sequence as in~\eqref{eq:identifiedexecution}.
Write $\widehat x_j=\widehat x_j^{h_i}(\bx(i),\bu_i)$
and $e_j=e_{j,t_i}^{h_i}(\bx(i),\bu_i)$.
For $0\le j<M$ with $x(iM+j)\in\mathcal D$,
equations~\eqref{eq:truestateL} and~\eqref{eq:retainederrorvalidity}
give $\norm{x(iM+j+1)-\widehat x_{j+1}}\le L_g\norm{x(iM+j)-\widehat x_j} + d_{t_i}^{\mathcal H}((\widehat x_j,u_{i,j}),h_i)$, where $\bu_i=(u_{i,0},\ldots,u_{i,M-1})$.
At $j=0$, the two states coincide.
If their distance at index $j$ is at most $e_j$, then
\eqref{eq:retainedpredictionerrors} bounds their distance at
index $j+1$ by $e_{j+1}$. Since
$\bu_i\in\mathcal S_{M,t_i}^{\mathcal H}(\bx(i),h_i)$,
equation~\eqref{eq:retainedsafeadmissibility} implies $x(iM+j+1)\in\Bball(\widehat x_{j+1},e_{j+1})
\subseteq\mathcal D$.
Induction, together with~\eqref{eq:modelupdateconditions}, proves $\norm{x(iM+j)-\widehat x_j}\le e_j\le\bar\Delta_j$, $x(iM+j)\in\mathcal D\qquad(0\le j\le M)$.
Therefore, we have $V(\bx(i+1)) \le V(\widehat x_M^{h_i}(\bx(i),\bu_i))+L_V\bar\Delta_M \le\rho(V(\bx(i)))+\Delta_0\le\rho(c)+\Delta_0<c$.

Starting from $\bx(0)\in\Omegac$, it follows by induction on $i$ that $x(k)\in\mathcal D$ for all $k\in\N$ and $\bx(i)\in\Omegac$ for all $i\in\N$.
The constants in~\eqref{eq:practicalratechoice} and~\eqref{eq:practicalthreshold} are independent of $i$ and $h_i$.
For $a\le V(\bx(i))\le c$, equation~\eqref{eq:thresholdgap} yields
$V(\bx(i+1)) \le\rho(V(\bx(i)))+\Delta_0 \le\rho(V(\bx(i)))+\vartheta\big(V(\bx(i))-\rho(V(\bx(i)))\big) = \tilde\rho(V(\bx(i)))$.
For $0\le V(\bx(i))\le a$, monotonicity of $\rho$ and~\eqref{eq:thresholdgap} give
$V(\bx(i+1)) \le\rho(a)+\Delta_0 = a-\frac{1-\vartheta}{\vartheta}\Delta_0 \le a$.
Combining these two inequalities yields $V(\bx(i+1))\le \max\{\tilde\rho(V(\bx(i))),a\}$ for all $i\in\N$.

At $i=0$, inequality~\eqref{eq:practicalcomparison} holds since $\bx(0)=\xi$.
Suppose that it holds at $i$. 
By monotonicity of $\tilde\rho$ and $\tilde\rho(a)\le a$, we obtain
$V(\bx(i+1)) \le\max\{\tilde\rho(V(\bx(i))),a\} \le\max\{\tilde\rho^{i+1}(V(\xi)),\tilde\rho(a),a\} = \max\{\tilde\rho^{i+1}(V(\xi)),a\}$.
This proves~\eqref{eq:practicalcomparison} by induction.
Finally, we get $\norm{x(iM+j)}\le b_j(\norm{\bx(i)})+\bar\Delta_j$, for $0\le j<M$.
Substitution into~\eqref{eq:betablock}--\eqref{eq:betainterpolation}
proves~\eqref{eq:intermiateTrajEstimate} and~\eqref{eq:practicalAsympticEstimate} with~\eqref{eq:practicalradius}.
\end{proof}

\begin{controlalgorithm}{Constraint-preserving data collection and model updates}
\label{alg:modelupdate}
\textit{Input Data:} $\Xset,\Uset,\mathcal D,V,M,\rho,b_j, L_g,L_V,c,\bar d$ satisfying the assumptions of
Proposition~\ref{prop:balltoset}, $r,\gamma$ from \eqref{eq:stagecost}--\eqref{eq:Jpi},
$\Vhat\in C(\Omegac)$, an initial data index $t_0\in\N$, an initial model $\widehat g=h_0$ satisfying
\eqref{eq:modelupdateconditions} at $t=t_0$,
and $x(0)\in\Omegac$.
Assume that the transition data \(\eqref{eq:transitiondata}\), the Lipschitz bound \(\eqref{eq:jointplantL}\), and a valid initial uncertainty set \(\mathcal H_0\) are available. Construct the uncertainty sets by \(\eqref{eq:nesteduncertainty}\).
Initialize $k=0$.
\begin{enumerate}[label=\arabic*.,leftmargin=1.8em]
\item Set $t=t_0+kM$ and measure $\xi_k=x(kM)$.
Compute $\mathcal H_t$ by \eqref{eq:nesteduncertainty}.
Fit a model $h$ to the available transition data.
\item Check Assumption~\ref{ass:standing_learned_syst}
for $h$ with the fixed $V,M,\rho,b_j$, and
\eqref{eq:modelupdateconditions} at $(t,h)$.
If both hold, set $\widehat g=h$.
Otherwise retain the previous $\widehat g$.
Set $d(z)=d^{\mathcal H_t}(z,\widehat g)$ by
\eqref{eq:retainedmodelerror}.
Form the bounds~\eqref{eq:retainedpredictionerrors}
and the safe set~\eqref{eq:retainedsafeadmissibility}
for the retained model.

\item Run Algorithm~\ref{alg:robustcontrol} at $\xi_k$
with $\widehat g,d,\Vhat$ to compute
$\bu^{\mathrm f}\in\mathcal S_{M,t}^{\mathcal H}(\xi_k,\widehat g)$.
If it returns failure, terminate without applying a new
control sequence.
Otherwise apply $u_j^{\mathrm f}$ at time $kM+j$,
$j=0,\ldots,M-1$, and record $z_{t+j+1} = (x(kM+j)^\top,(u_j^{\mathrm f})^\top)^\top$, $y_{t+j+1}=x(kM+j+1)+w_{t+j+1}$.
\item Append
$(z_{t+j+1},y_{t+j+1},\nu_{t+j+1})$ with
$\nu_{t+j+1}\ge\norm{w_{t+j+1}}$ to
\eqref{eq:transitiondata}, for $j=0,\ldots,M-1$.
Set $k\leftarrow k+1$ and return to step 1.
\end{enumerate}
\end{controlalgorithm}

\begin{figure}
    \centering
    \includegraphics[width=1\linewidth]{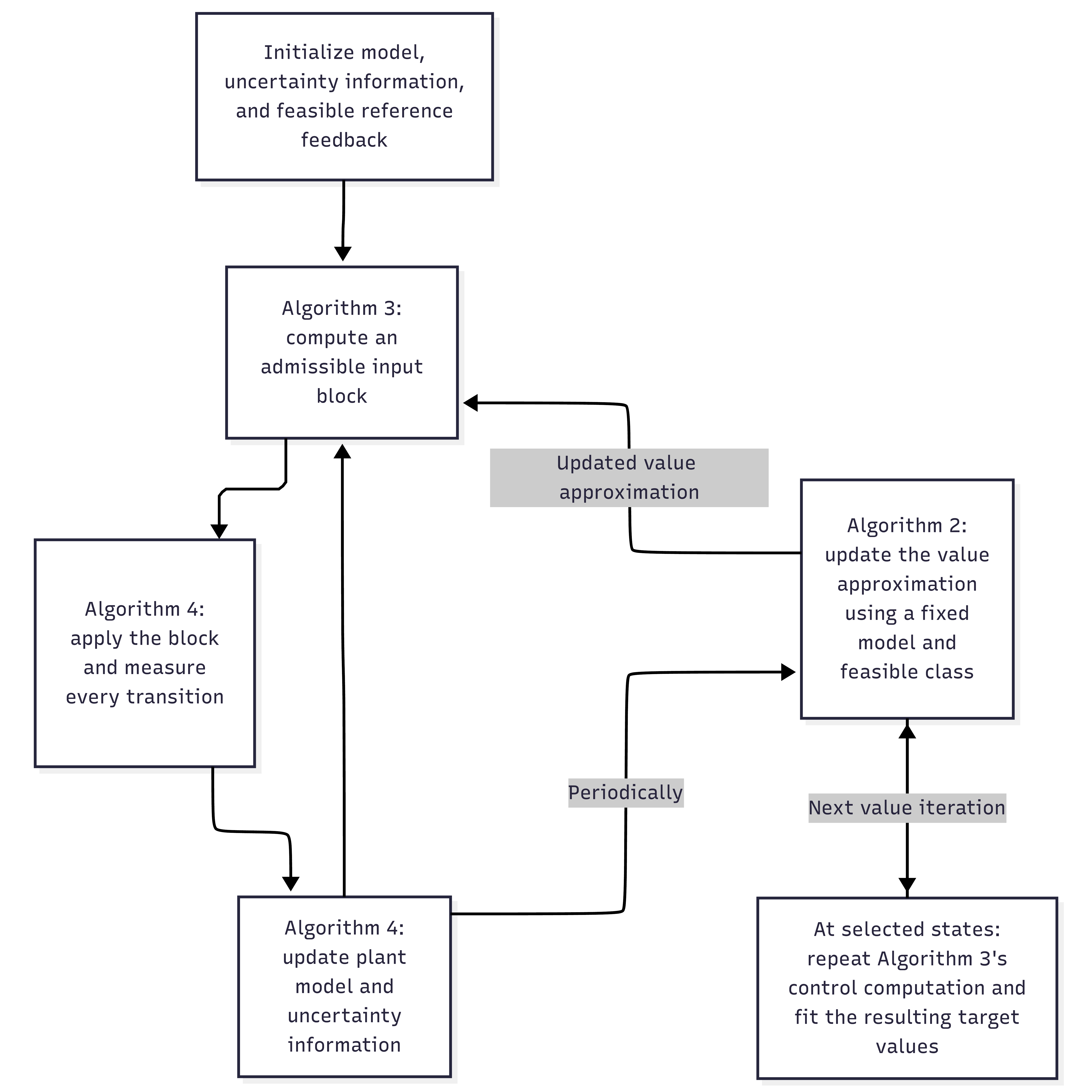}
    \caption{Computation of learned finite-step control laws through Algorithms~\ref{alg:valueiteration} to~\ref{alg:modelupdate}.}
    \label{fig:workflow}
\end{figure}

Figure~\ref{fig:workflow} summarizes the workflow of our learned fs-CLF control design through Algorithms~\ref{alg:valueiteration} to~\ref{alg:modelupdate}.

\subsection{Verification of the assumptions for the two benchmarks}
\label{sec:additionalbenchmarkchecks}
We verify the assumptions of Theorem~\ref{thm:practical} and
Algorithms~\ref{alg:robustcontrol}--\ref{alg:modelupdate} before
computing the numerical results. The verification of
\eqref{eq:safenonempty} below concerns every state in $\Omegac$.
It is independent of the choice of the continuous approximation
$\Vhat$ and of the initial conditions used in the numerical results.

\subsubsection{Single integrator with parameter uncertainty}
\label{sec:singleintegratorchecks}
Let $x=(x_1,x_2)$, $u=(u_1,u_2)$, and
$\dot x=(1+\lambda)u$, where the unknown parameter
$\lambda\in(-0.1,0.1)$ is constant.
With a constant input over each sampling interval,
\begin{equation}\label{eq:siprecheckmodel}
g(x,u)=x+T(1+\lambda)u,\qquad
\widehat g(x,u)=x+T\widehat a u,
\end{equation}
where $\widehat a=1+\widehat\lambda$.

The desired trajectory $x_d(t)$ is given to the tracking controller.
BARN supplies the geometric waypoints~\cite{Perille2020BARN}. Let $P:[0,L]\to\R^2$ denote their polygonal interpolation in metres, including the segments to the prescribed start and goal, parametrized by distance from the start.
For this experiment we prescribe a speed $v_d=0.55$ m/s along the path until the goal and zero speed thereafter.
The distance travelled by time $t$ is $\min\{v_dt,L\}$, so the given reference is $x_d(t) := P\big(\min\{v_dt,L\}\big)$, $t\ge0$.
At sample $k$ the remaining distance is $\ell_{\rm r}(k):=L-\min\{v_d kT,L\}$.
In that way, $\ell_{\rm r}(0)=L$.
Hence $p(\ell_{\rm r})=P(L-\ell_{\rm r})$ for $0\le\ell_{\rm r}\le L$, and $p(\ell_{\rm r}(k))=x_d(kT)$.
Extend $p$ by its endpoint values outside $[0,L]$.
Letting $e(k):=x(k)-x_d(kT)$, subtracting consecutive reference positions in~\eqref{eq:siprecheckmodel} gives $e(k+1)=e(k)+Ta u(k)-[x_d((k+1)T)-x_d(kT)]$.
Define $z=(y,\ell_{\rm r})$, where $y:=x-x_{\rm G}$, $x_{\rm G}:=P(L)$.
Over one sample the reference advances by $h=Tv_d$ or by the distance left to the goal, whichever is smaller.
Together with $y^+=y+Ta u$, this gives
\begin{equation}\label{eq:barnsystem}
\begin{aligned}
 F(\ell_{\rm r})&=\max\{0,\ell_{\rm r}-h\},\qquad \ell_{\rm r}\ge0,\\
 z^+=G_a(z,u)&=(y+Ta u,F(\ell_{\rm r})),\\
 e&=y+x_{\rm G}-p(\ell_{\rm r}).
\end{aligned}
\end{equation}In the experiment $T=0.20$ s.
Let $a_\star=1+\lambda_\star$ and $\widehat a$ belong to $[a_-,a_+]=[1.039,1.041]$. For $M=4$ we use
\begin{equation}\label{eq:barndata}
\begin{gathered}
 d_l=(\cos(l\pi/4),\sin(l\pi/4)),\quad N(e)=\max_{0\le l<8}d_l^\top e,\\
 V(z)=c\big(N(e)/R+|\ell_{\rm r}|/S\big),\\
 c=10,\quad R=0.10\ {\rm m},\quad S=15.00\ {\rm m},\\
 \Omegac=\{V\le10\},\qquad \mathcal D=\Xset=\{V\le11\},\\
 \Uset=[-u_{\max},u_{\max}]^2,\quad u_{\max}=0.68\ {\rm m/s},\\
 q=1-h/(1.10S),\qquad \rho(t)=q^Mt,\qquad b_j(t)=4t,\\
 r(u)=0.60\norm{u}_1,\qquad\gamma=0.95 .
\end{gathered}
\end{equation}
The term $|\ell_{\rm r}|/S$ makes $V$ positive even at zero tracking error until the reference reaches the goal. 
Both $p$ and $F$ are $1$-Lipschitz, and $\cos(\pi/8)\norm{e}\le N(e)\le\norm{e}$. Thus $\norm{z}\le\norm{e}+2|\ell_{\rm r}|\le(2S/c)V(z)$. By the reverse triangle inequality for $N$, a Lipschitz constant for $V$ is $L_V=c\sqrt{R^{-2}+(R^{-1}+S^{-1})^2}$. This gives~\eqref{eq:properV} with $\alpha_1(t)=ct/(2S)$ and $\alpha_2(t)=L_Vt$ and implies compactness of the two sublevels. From~\eqref{eq:barnsystem}, $G_a(0,0)=0$, $L_g=1$, and $L_{\rm sm}=\sqrt{1+(Ta_+)^2}$ is a joint Lipschitz bound. In~\eqref{eq:K_bounded_g_open_loop} we may take $\kappa_1(t)=t$ and $\kappa_2(t)=Ta_+t$.

Fix $z\in\mathcal D$ and $a\in[a_-,a_+]$. To prove nonemptiness, choose $Ta u=p(F(\ell_{\rm r}))-p(\ell_{\rm r})-\eta e$, where $\eta=0.15$. Then $e^+=(1-\eta)e$. 
With $r_D=1.10R/\cos(\pi/8)$, we have $Ta\norm{u}_\infty\le h+\eta r_D<0.128\ {\rm m}<Ta_-u_{\max}$.
Since $|\ell_{\rm r}|\le1.10S$, $|F(\ell_{\rm r})|\le q|\ell_{\rm r}|$ and $1-\eta\le q$, this input gives $V(z^+)\le qV(z)$. Repetition gives $V(z_j)\le q^jV(z_0)$ and keeps all successors in $\mathcal D$. Moreover, $\norm{y_j}\le\norm{e_0}+|\ell_{\rm r}(0)|\le\norm{y_0}+2|\ell_{\rm r}(0)|$ and $|\ell_{\rm r}(j)|\le|\ell_{\rm r}(0)|$, so $\norm{z_j}\le4\norm{z_0}$. This proves admissibility for the true and fitted models in~\eqref{eq:admset} and~\eqref{eq:practicaladmissible}.

For Algorithm~\ref{alg:modelupdate}, the prior segment $\{G_a(z,u):a\in[a_-,a_+]\}$ and each measurement ball in~\eqref{eq:nesteduncertainty} contain the true successor. For both the fixed and updated models the error is therefore bounded by $T(a_+-a_-)\norm{u}\le\bar d:=T(a_+-a_-)\sqrt{2}\,u_{\max}$. Since $L_g=1$, $\bar\Delta_M=M\bar d<1.54$ mm and $\Delta_0=L_V\bar\Delta_M\simeq0.21833<c-\rho(c)\simeq0.26401$.
For $z_0\in\Omegac$ the constructed predictions satisfy $V(z_j)+L_V\bar\Delta_j\le c+\Delta_0<11$ for $1\le j\le M$. Hence their prediction balls lie in $\mathcal D$, proving~\eqref{eq:safenonempty} and~\eqref{eq:modelupdateconditions}. Theorem~\ref{thm:practical} gives practical asymptotic stability of the augmented origin with threshold $a\simeq9.05$ in~\eqref{eq:practicalthreshold}. Since $\norm{e}\le R V(z)/(c\cos(\pi/8))$, the nominal tracking error converges to zero. For the uncertain plant,~\eqref{eq:practicalcomparison} gives $\limsup_{i\to\infty}\norm{e(iM)}\le Ra/(c\cos(\pi/8))<0.099$ m.

\subsubsection{Torque-controlled inverted pendulum}
Let $x=(\theta,\omega)$ and
\begin{equation}\label{eq:pendprecheckcontinuous}
\dot\theta=\omega,\qquad
\dot\omega=A\sin\theta-B\omega+Cu,
\end{equation}
where the uniform rod has mass $m=1.00$ kg, length $l=1.00$ m,
and $g_{\rm grav}=9.81\ {\rm m/s^2}$. Its moment of inertia about
the pivot is $ml^2/3$, and its center of mass is at $l/2$. Hence $A_\star=\frac{3g_{\rm grav}}{2l}$, $B=0$, $C=\frac{3}{ml^2}=3.00$, $\Uset=[-5.50,5.50]$.
The torque $u$ is measured in N\,m.

In~\eqref{eq:pendprecheckcontinuous} take $A=A_\star$ for the true plant.
Define $S(v)=\min\{8.00,\max\{-8.00,v\}\}$. 
Given the sampling period $T=0.05$ s, the true sampled dynamics and the fitted family are
\begin{equation}\label{eq:g2plant}
\begin{aligned}
v_A(x,u)&=S\big(\omega+T(A\sin\theta+3u)\big),\\
h_A(x,u)&=(\theta+Tv_A(x,u),\,v_A(x,u)),\\
g&=h_{A_\star},\quad A_\star=\frac{3g_{\rm grav}}{2l},\\
A&\in I=[14.69,14.74].
\end{aligned}
\end{equation}
Pick the discrete-time true dynamics $g$ as $g=h_{A_\star}$ in~\eqref{eq:g2plant}.
In this update, the angular speed is updated first, clipped to $[-8.00,8.00]$ rad/s, and then used in the angle update. The fitted map is $\widehat g=h_{\widehat A}$ with $\widehat A\in I$ with $I$ as in~\eqref{eq:g2plant}.

Choose
\begin{equation}\label{eq:g2clf}
\begin{aligned}
P&=\begin{bmatrix}2.00&2.50\\2.50&6.25\end{bmatrix},\qquad
 V(x)=x^\top P x,\\
\mathcal D=\Xset&=\{V\le17.00\},\qquad \Omega_{16}=\{V\le16.00\},\\
\rho(s)&=0.98^M s,\qquad b_j(s)=3.00s\quad(1\le j\le M).
\end{aligned}
\end{equation}
Equivalently, $V(x)=\theta^2+(\theta+2.50\omega)^2$.
In that way, condition~\eqref{eq:properV} holds on $\Xset$ with $0.84I_2\prec P\prec7.41I_2$.

\emph{Nonemptiness of the admissible sets.}
We verify~\eqref{eq:admset} and~\eqref{eq:practicaladmissible} by considering the whole interval of attainable next states. Put $\kappa=0.40$, $s=\theta+\omega/\kappa$, $d=1+T\kappa$, and $\lambda=0.98$. As $u$ ranges over $\Uset$, the next coordinate $s^+=\theta+(T+1/\kappa)v_A(x,u)$ ranges over a closed interval. Substituting $\omega=\kappa(s-\theta)$ shows that $\lambda s$ belongs to this interval if and only if
\begin{equation}\label{eq:g2interval}
\left|\frac{\kappa^2}{d}\theta+
 \frac{\kappa}{T}\left(\frac{\lambda}{d}-1\right)s
 -A\sin\theta\right|\le 3\max_{u\in\Uset}|u|.
\end{equation}
Because $\theta^2+s^2\le17.00$, the left-hand side is at most
\begin{equation}\label{eq:g2availability}
14.74+\sqrt{17}\left\|
\left(\frac{\kappa^2}{d},
 \frac{\kappa}{T}\left(\frac{\lambda}{d}-1\right)\right)
\right\|<16.20<16.50.
\end{equation}
In the coordinates $(\theta,s)$, $V=\theta^2+s^2$. For the particular attainable successor $s^+=\lambda s$, substituting $v_A=\kappa(s^+-\theta^+)$ into $\theta^+=\theta+Tv_A$ yields
\begin{equation}\label{eq:g2matrix}
 \binom{\theta^+}{s^+}=F\binom{\theta}{s},\qquad
 F=\begin{bmatrix}1/d&T\kappa\lambda/d\\0&\lambda\end{bmatrix}.
\end{equation}
We note tha $F^\top F\prec0.98I_2$.
For every $A\in I$, the attainable-input inequality~\eqref{eq:g2interval} follows from~\eqref{eq:g2availability} on this whole domain. Equation~\eqref{eq:g2matrix} therefore gives, at every $x\in\mathcal D$, an input in $\Uset$ for which
$V(x^+)\le0.98V(x)$. Repeating this argument $M$ times gives
$V(x_M)\le\rho(V(x))$ and
$\|x_j\|\le\sqrt{7.41/0.84}\,\|x\|\le b_j(\|x\|)$
for $1\le j\le M$.

For two state--input pairs, denote differences by $\delta$ and put $b=14.74T|\delta\theta|+|\delta\omega|+3T|\delta u|$. Since $S$ and $\sin$ are $1$-Lipschitz,~\eqref{eq:g2plant} gives $|\delta v_A|\le b$, $|\delta\theta^+|\le|\delta\theta|+Tb$.
The map from $(|\delta\theta|,|\delta\omega|,|\delta u|)$ to $(|\delta\theta|+Tb,b)$ has Euclidean operator norm below $1.48$ and its restriction to $\delta u=0$ has norm below $1.47$. Thus $L_g=1.47$, $L_{\rm sm}=1.48$, and~\eqref{eq:K_bounded_g_open_loop} holds with $\kappa_1(t)=1.47t$, $\kappa_2(t)=3T\sqrt{1+T^2}\,t$. On the convex set $\mathcal D$, $\norm{2Px}\le2\sqrt{17(7.41)}<22.45=L_V$.

At $z=(x,u)$, the prior $\mathcal H_0(z)$ is the sum of the segment joining $h_{14.69}(z)$ to $h_{14.74}(z)$ and the segment $\{t v:|t|\le10^{-3}\}$, where $v=(T,1)/\sqrt{1+T^2}$. It contains $g(z)$ since $v_A(x,u)$ is monotone in $A$ for fixed $z$. For every $A\in I$ and every $y\in\mathcal H_0(z)$, we get $\norm{y-h_A(z)}\le0.05T\sqrt{1+T^2}+10^{-3}<\bar d=0.0036$.
Each observation ball in~\eqref{eq:nesteduncertainty} contains $g(z)$ by~\eqref{eq:transitiondata} and the joint Lipschitz bound.
Its intersection with the prior therefore preserves the true successor and the latter bound $\norm{y-h_A(z)}$ at every update.

For $\xi\in\Omegac$, the existence sequences above satisfy
$V(\widehat x_j)\le16.00$ for all $1\le j\le M$. By the
triangle inequality in the $P$-norm, every
$y\in\Bball(\widehat x_j,e_j)$ satisfies
$\sqrt{V(y)}\le4+\sqrt{7.41}\,\bar\Delta_M$.
Thus~\eqref{eq:safenonempty} and~\eqref{eq:uniformmismatchsmallness}
follow from the two checks
\begin{equation}\label{eq:pendprechecksafeconstraints}
 \bar\Delta_M<\frac{\sqrt{17}-4}{\sqrt{7.41}},\qquad
 L_V\bar\Delta_M<16(1-0.98^M),
\end{equation}
where $\bar\Delta_M=\bar d\sum_{j=0}^{M-1}1.47^j$ by
\eqref{eq:uniformpredictionerrors}.
Already at $M=1$, $\bar\Delta_1=0.0036<(\sqrt{17}-4)/\sqrt{7.41}$ and $L_V\bar\Delta_1=0.081<0.32=16(1-0.98)$. 
Table~\ref{tab:pendcertificate} gives upper bounds on $10^3\bar\Delta_M$ and $\Delta_0=L_V\bar\Delta_M$, and lower bounds on $c-\rho(c)$. The first column of bounds is less than
$45.00$ in every row, while $10^3(\sqrt{17}-4)/\sqrt{7.41}>45$.
\begin{table}[t]
\centering
\caption{Pendulum certificates at $T=0.05$ s. Error bounds are
rounded upwards and contraction budgets downwards.}
\label{tab:pendcertificate}
\begin{tabular}{rrrr}\hline
$M$&$10^3\bar\Delta_M$&$\Delta_0$&$c-\rho(c)$\\\hline
2&8.90&0.20&0.63\\
3&16.68&0.38&0.94\\
4&28.11&0.64&1.24\\
5&44.92&1.01&1.53\\\hline
\end{tabular}
\end{table}
Both inequalities in~\eqref{eq:pendprechecksafeconstraints}
hold uniformly over the retained model family. This establishes
both parts of~\eqref{eq:modelupdateconditions} at every update.
For $M=4$, Theorem~\ref{thm:practical} gives
$x(k)\in\mathcal D$ at all sampling instants and
$x(iM)\in\Omegac$ at every block boundary under
Algorithms~\ref{alg:robustcontrol} and~\ref{alg:modelupdate}.
In the experiment, the threshold in~\eqref{eq:practicalthreshold} is $a\simeq10.78$.

\section{Numerical Results: Planar Trajectory Tracking}
\label{sec:barntracking}

We consider trajectory tracking for the single integrator~\eqref{eq:siprecheckmodel}
with an uncertain input coefficient. The reference trajectories are prescribed
in worlds~5, 34 and~53 of the Benchmark for Autonomous Robot Navigation
(BARN)~\cite{Perille2020BARN}. Each experiment uses the
released  obstacle positions, cylinder radii, start and target
positions, and geometric path. 
We use the static BARN paths with a prescribed timing common to
all three experiments. 

In the world frame the start and goal are $x_{\rm S}=(-2.25,3.00)$ m and $x_{\rm G}=(-2.25,13)$ m.
The robot disk of radius $r_{\rm B}=334$ mm encloses the
$0.508\times0.43$ m rectangular footprint reported
in~\cite{Perille2020BARN}. All obstacles are vertical cylinders of
radius $75$ mm. We use the observation rectangle
$[-4.6,0.1]\times[0,14]$ m, which contains the complete runs.
The ground plane extends beyond this rectangle.
The true plant has $a_\star=1+\lambda_\star=1.0408$.
In a separate initial measurement at $z=0$, a command
$u=(0.60,0)$ m/s gives an observed first-coordinate displacement of
$124.80$ mm with Euclidean output error bounded by $0.12$ mm.
Intersecting the prior $[0.90,1.10]$ with this measurement gives
$[a_-,a_+]=[1.039,1.041]$ and $\widehat a_0=1.040$.
All runs start with $x(0)-x_{\rm S}=(12,-15)$ mm. The three initial values of $V$ are
$9.82$, $9.22$ and $9.08$, respectively.

For each world, Algorithm~\ref{alg:valueiteration} starts from
$\Vhat_0=0$ and uses the origin and 31 further samples from
$\Omegac$, with equal weights and four iterations.
We fit $W_\theta(z)=\theta_1 |e_1|/R+\theta_2 |e_2|/R+\theta_3 |\ell_{\rm r}|/S$, $\theta\ge 0$,
by nonnegative least squares in~\eqref{eq:generalfit}, followed by
\eqref{eq:projection}. In~\eqref{eq:costbound} we may use
$R_c=2S$ and $\bar r=1.20u_{\max}$. On $\Omegac$,
$W_\theta\le U(\theta):=\max\{\theta_1+\theta_2,\theta_3\}$ since
$|e_i|\le N(e)$. The computed $U(\theta)$ lies below $\bar V_c$
at every iteration. The transition model in the offline
Problem~\ref{prob:approximate_optimal_control} is the true simulator.
For $U_i=U(\theta_i)$ put $Q_i^+=M(11+\bar r)+\gamma U_i$.
The all-state feasible block established in
Section~\ref{sec:singleintegratorchecks} and nonnegativity of the
cost give the uniform bounds $\bar\varepsilon_i=\max\{U_{i+1},Q_i^+\}$,
$\bar\delta_i=Q_i^+$ and $\bar R_{i+1}=\max\{U_{i+1},Q_{i+1}^+\}$
in~\eqref{eq:verifiedupdate}--\eqref{eq:verifiedresidual}.
The last residual bounds are $187.23$, $181.92$ and $179.46$.
The iteration limit determines the stopping point.
The final parameter vectors in world order are $(4.78,3.25,147.34)$, $(5.88,4.79,141.74)$ and $(5.10,6.12,139.15)$.
For each world let $\mathcal Z$ denote its set of $32$ sampled augmented states.
The final mean squared losses $|\mathcal Z|^{-1}\sum_{\zeta\in\mathcal Z}|\Vhat_4(\zeta)-q_4(\zeta)|^2$, where $q_4$ denotes the target~\eqref{eq:targets} computed with $\Vhat_3$, are $0.8$, $1.27$ and $1$.

\begin{figure*}[!t]\centering
\includegraphics[width=0.74\textwidth]{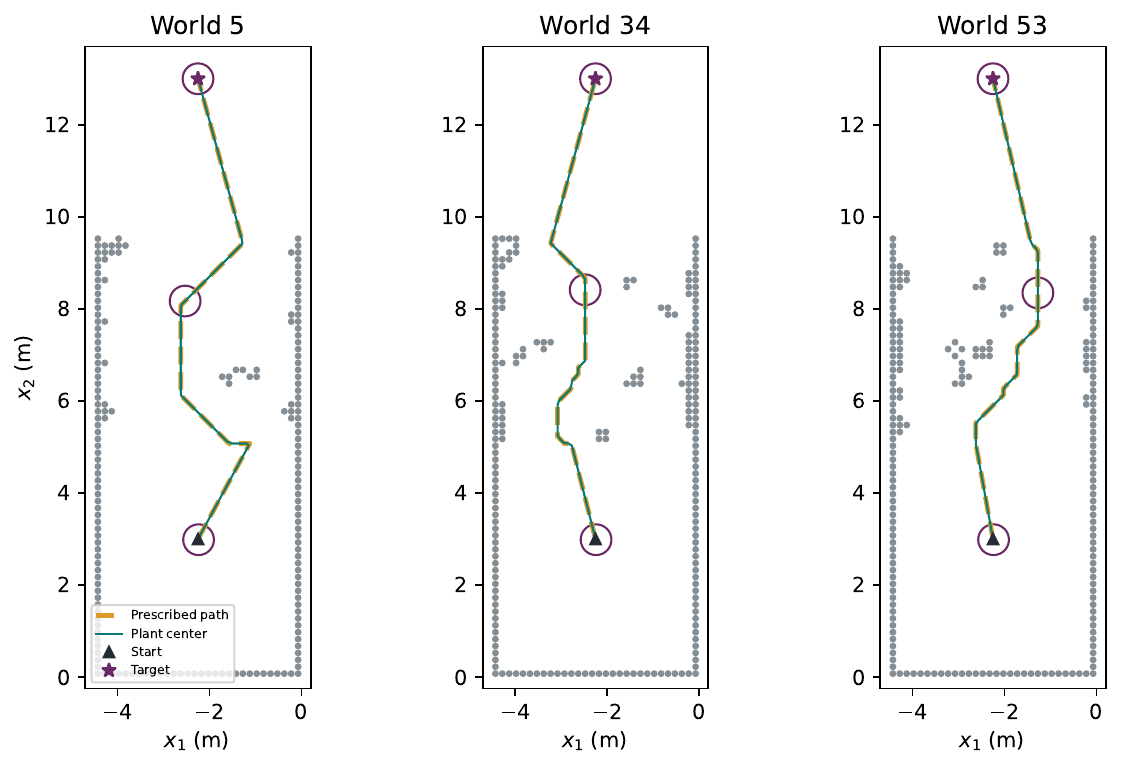}
\caption{Algorithm~\ref{alg:modelupdate} on the three static BARN paths. Circles give the enclosing unicycle disk at the initial, middle and final samples.}\label{fig:barnpaths}
\end{figure*}

\begin{figure*}[!t]\centering
\includegraphics[width=0.86\textwidth]{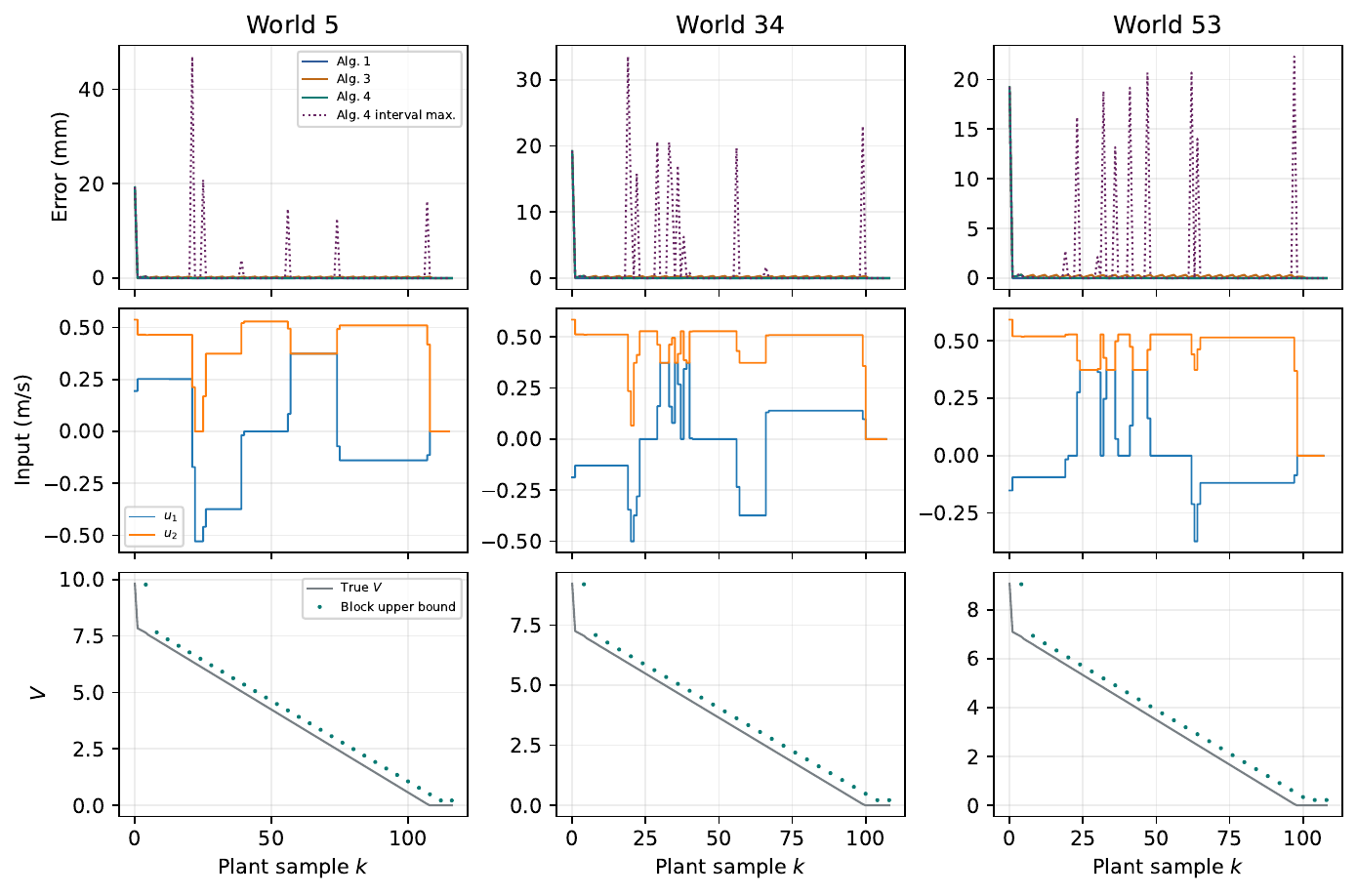}
\caption{Tracking errors for the three controllers, both applied velocity components for Algorithm~\ref{alg:modelupdate}, and $V$ along the true trajectory. The dotted error trace gives the maximum over each held-input interval. The points in the last row give $\rho(V(z(iM)))+\Delta_0$ at the following block endpoint.}\label{fig:barntrack}
\end{figure*}

\begin{figure*}[!t]\centering
\includegraphics[width=0.82\textwidth]{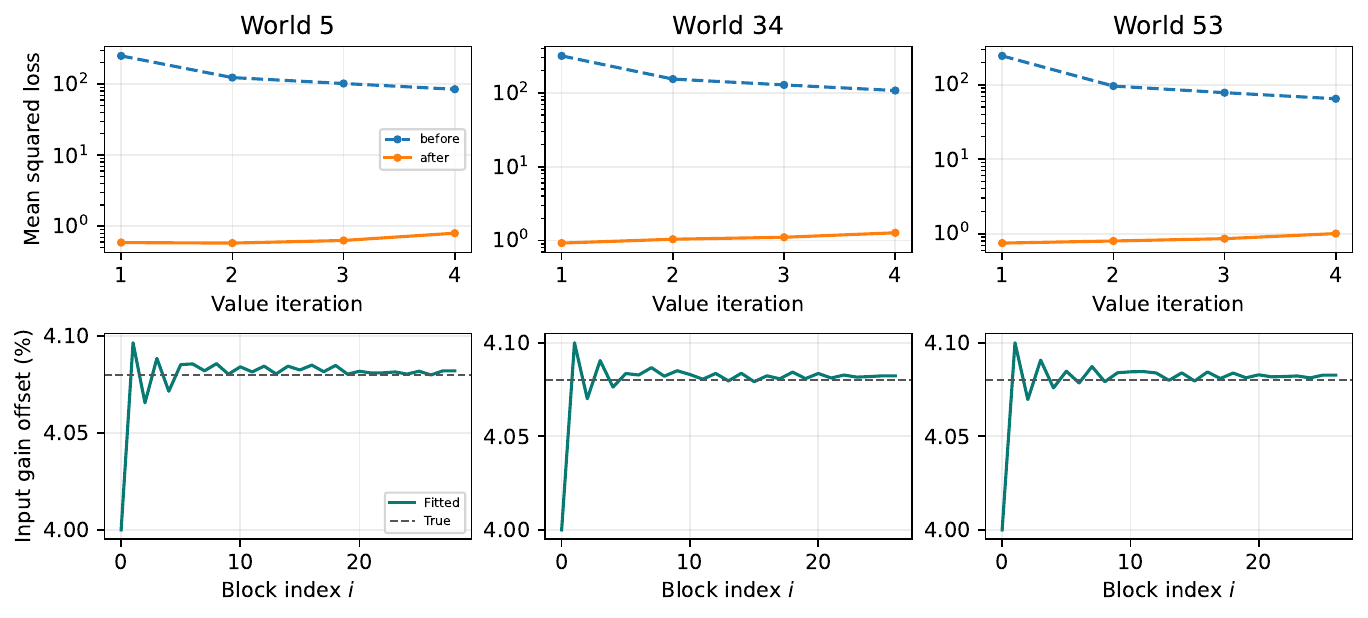}
\caption{Sample losses on the current targets before and after each fit, and the retained offset $100\widehat\lambda_i$ in percent under Algorithm~\ref{alg:modelupdate}.}\label{fig:barnlearn}
\end{figure*}

Table~\ref{tab:barnruns} reports the executed cost $J^{\rm ap}(NM)=\sum_{i=0}^{N-1}\gamma^i\sum_{j=0}^{M-1}\ell(z(iM+j),u(iM+j))$ for $N$ complete blocks. At each sample $k$
we evaluate $V(z(k))+r(u(k))$ from~\eqref{eq:barndata} at
the true augmented state and the applied input. 
The table also gives the sampled root-mean-square position error including the initial sample, the largest continuous-time tracking error, and the minimum separation between the robot disk and the obstacles along the true trajectory.
Following a zero-order-hold implementation, the input is held constant over each sampling interval.
In that way, the minimum separation is computed along the line segments connecting successive positions of the true plant.
Algorithm~\ref{alg:modelupdate} gives maximum center differences $\max_k\norm{x^{\rm U}(k)-x(k)}$ of $0.34$, $0.35$ and $0.35$ mm.
Since both centers are affine during each held input, the maximum center difference on an interval is attained at an endpoint. Even the largest difference across all nine runs ($7.7$ mm) is below the smallest certified separation ($12.81$ mm).
Figures~\ref{fig:barnpaths} to~\ref{fig:barnlearn} show the paths, the sample-by-sample control and the identification results.

\begin{table*}[!t]\centering
\caption{BARN trajectory tracking. The last three columns give mean/maximum elapsed times per block in milliseconds. First command includes model fitting; complete update also includes plant execution and data collection on the simulation host.}
\label{tab:barnruns}
\begin{tabular}{rrrrrrrrrr}\hline
World & Alg. & Samples applied & $J^{\rm ap}$ & RMS (mm) & Peak (cm) & Min. separation (cm) & Inner solve & First command & Complete update\\\hline
5 & 1 & 116 & 317.91 & 1.78 & 4.69 & 49.10 & 3.83/4.62 & 3.83/4.62 & 4.06/4.81\\
5 & 3 & 116 & 319.10 & 1.79 & 4.69 & 49.10 & 3.93/4.28 & 3.93/4.28 & 4.15/4.54\\
5 & 4 & 116 & 318.07 & 1.78 & 4.69 & 49.10 & 4.24/5.00 & 4.29/5.06 & 4.53/5.33\\
34 & 1 & 108 & 280.68 & 1.84 & 3.35 & 19.72 & 3.59/3.94 & 3.60/3.94 & 3.79/4.15\\
34 & 3 & 108 & 281.83 & 1.85 & 3.34 & 19.72 & 3.72/4.22 & 3.72/4.22 & 3.92/4.40\\
34 & 4 & 108 & 280.83 & 1.84 & 3.35 & 19.72 & 4.19/4.57 & 4.24/4.59 & 4.46/4.88\\
53 & 1 & 108 & 271.87 & 1.84 & 2.23 & 18.69 & 3.85/7.73 & 3.85/7.73 & 4.05/7.94\\
53 & 3 & 108 & 273.04 & 1.85 & 2.22 & 18.69 & 3.69/3.97 & 3.69/3.97 & 3.89/4.18\\
53 & 4 & 108 & 272.03 & 1.84 & 2.23 & 18.69 & 4.19/4.48 & 4.24/4.53 & 4.47/4.97\\
\hline\end{tabular}\end{table*}

\section{Numerical Results: Large-Angle Swing-Up}
\label{sec:gymswing}
We use the sampled dynamics of the \href{https://gymnasium.farama.org/environments/classic_control/pendulum/}{Gymnasium Pendulum-v1 benchmark}~\cite{GymnasiumBenchmark}.
In Gymnasium version 1.2.2 we set both the declared action space and the internal torque clipping limit to the selected rating. 
The measured state is the simulator's internal $x=(\theta,\omega)$, with the angle in radians relative to the upright position and speed in rad/s. 
We use these numerical coordinates in the Euclidean state norm. 
We report simulation results for two test following initial conditions
\begin{equation}\label{eq:g2starts}
\xi_1=(2.75,0),\qquad \xi_2=(-1.6,0.65).
\end{equation}
Hence, their initial angles are $157.56^\circ$ and $-91.67^\circ$.
For the chosen initial conditions, we obtain $V(\xi_1)\simeq15.13$ and $V(\xi_2)\simeq2.56$.

\begin{figure}[!htbp]\centering
\includegraphics[width=\linewidth]{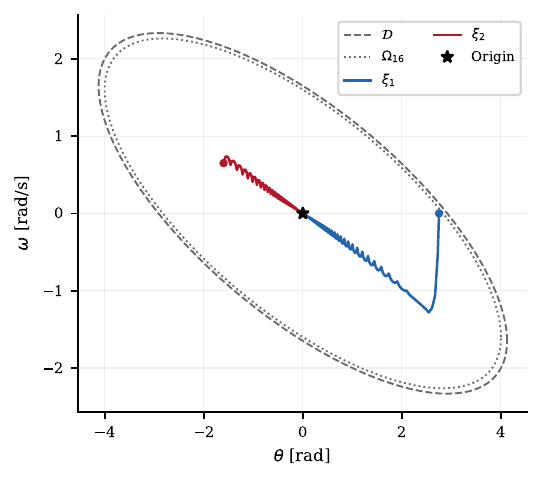}
\caption{Algorithm~\ref{alg:modelupdate}, polynomial value approximation. State trajectories from~\eqref{eq:g2starts}, the model domain $\mathcal D$ and the block-boundary set $\Omega_{16}$ in~\eqref{eq:g2clf}.}
\label{fig:g2sets}\end{figure}

\begin{figure*}[!htbp]\centering
\includegraphics[width=\linewidth]{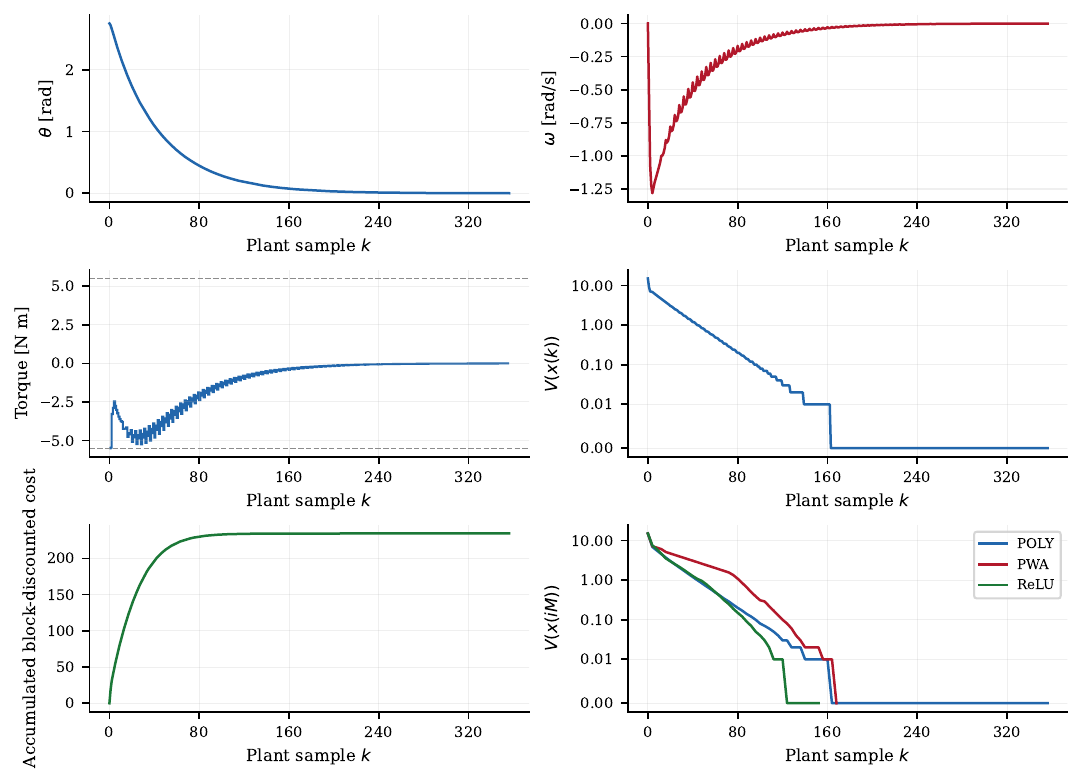}
\caption{Pendulum from $\xi_1$, $T=0.05$ s, $M=4$, $\gamma=0.95$. Algorithm~\ref{alg:modelupdate} with the polynomial fit applies all 356 recorded torques.}
\label{fig:g2seq}\end{figure*}

We collect 64 one-step transitions by independently drawing states uniformly in $\{V\le4\}$ and torques uniformly in $[-1,1]$ N\,m. Each prediction ball is checked to lie in $\mathcal D$ before applying its torque. The offline speed measurement includes an offset $\sin(\theta)$ mrad/s. In the numerical state coordinates the observation bounds~\eqref{eq:transitiondata} are $10^3\nu_s=2$ offline and $10^3\nu_s=1$ online. 
Least squares fits $\omega^+-\omega-3Tu$ to $TA\sin\theta$, with $A\in I$.
The initial fit is $\widehat A_0\simeq14.74$. On the training observations the angle and speed RMSEs are $0.04$ mrad and $0$ mrad/s to displayed precision. On 32 separate observations without added bias they are $0.03$ mrad and $0.69$ mrad/s.

We initialize $\Vhat_0=0$. Define the normalized state
$z=(\theta/\sqrt{17},\omega/(0.40\sqrt{34}))$.
For each direction
$d\in\{(1,0),(0,1),(1,1)/\sqrt{2},(1,-1)/\sqrt{2}\}$
the polynomial family uses the feature $(d^\top z)^2$.
For each $a\in\{0,0.25,0.5\}$ the continuous PWA family uses
$\max\{|d^\top z|-a,0\}$. The ReLU family uses the two features
$\max\{\sigma d^\top z-a,0\}$, $\sigma\in\{-1,1\}$. For Algorithm~\ref{alg:valueiteration}, use one sample at the origin and 23 independent samples drawn uniformly in $\Omega_{16}$, with equal weights and four iterations.
In this benchmark the simulator supplies the true map $g$ for the offline Problem~\ref{prob:approximate_optimal_control} and Algorithm~\ref{alg:control}. The online predictions in~\eqref{eq:modelcontrolobjective} use the fitted map. Each feasible nominal block gives a target $y^{(s)}$ by~\eqref{eq:targets}. Its first stage yields the sample gap bound $y^{(s)}-V(\xi^{(s)})$ in~\eqref{eq:sampletargetgap}. We solve~\eqref{eq:generalfit} by nonnegative least squares, then project by~\eqref{eq:projection}. In~\eqref{eq:costbound} take $R_c=\sqrt{c/0.84}$ and $\bar r=0.1(5.5)^2$, giving $\bar V_c\simeq76779$ in~\eqref{eq:constantupper}.

Let $\mathcal Z=\{\xi^{(s)}:1\le s\le24\}$ denote the fixed training set.
At fitting iteration $i\in\{1,2,3,4\}$, let $q_i(\xi^{(s)})$ be the
target~\eqref{eq:targets} computed with $\Vhat_{i-1}$.
For any value approximation $W$, the reported loss on these targets is
\begin{equation}\label{eq:newfitloss}
E_i(W)=\sum_{\zeta\in\mathcal Z}
                |W(\zeta)-q_i(\zeta)|^2,\qquad |\mathcal Z|=24.
\end{equation}
We evaluate~\eqref{eq:newfitloss} at $W=\Vhat_{i-1}$ before fitting
and at $W=\Vhat_i$ after fitting. Both evaluations use the same
$q_i$. The targets are recomputed at the next iteration.

For the uniform checks of Algorithm~\ref{alg:valueiteration}, let $U_i$ bound the projected fit above on $\Omega_{16}$. Use $U_0=0$, then take the minimum of $\bar V_c$ and the sum of the nonnegative coefficients times their respective feature bounds. The four absolute linear forms in the normalized coordinates are bounded by $1,1,\sqrt2,\sqrt2$ before applying the square or hinge. Therefore any feasible objective lies between zero and $Q_i^+=M(17.00+\bar r)+\gamma U_i$. Using the existence result~\eqref{eq:g2matrix} in~\eqref{eq:verifylower}--\eqref{eq:verifygap} yields $\bar\varepsilon_i=\max\{U_{i+1},Q_i^+\}$, $\bar\delta_i=Q_i^+$, and $\bar R_{i+1}=\max\{U_{i+1},Q_{i+1}^+\}$ in~\eqref{eq:verifiedresidual}. The fourth residual upper bounds are approximately $579.70$, $406.89$ and $722.35$. The stopping condition is the iteration limit. Table~\ref{tab:g2fit} gives the losses on the fourth set of targets. These sampled losses and the uniform bounds quantify different errors.
\begin{figure*}[!htbp]\centering
\includegraphics[width=\linewidth]{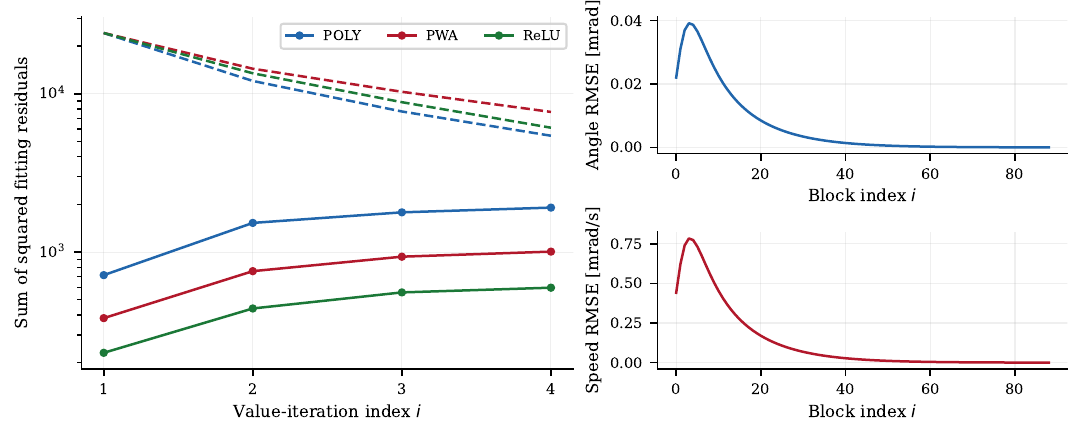}
\caption{Value-fitting losses~\eqref{eq:newfitloss}, dashed before and solid after fitting (left); angle and angular-speed residuals on each newly collected block for Algorithm~\ref{alg:modelupdate} with the polynomial fit from $\xi_1$ (right).}
\label{fig:g2learning}\end{figure*}

\begin{table}[!htbp]\centering
\caption{Value approximation. Fourth-iteration loss before and after fitting, total regression time over four iterations and maximum target-solve time (ms).}
\label{tab:g2fit}
\begin{tabular}{lrrrr}\hline
Family&Before&After&Regression&Target max.\\\hline
POLY & 5429.83 & 1909.76 & 0.11 & 21.74 \\
PWA & 7677.62 & 1006.81 & 0.15 & 80.40 \\
ReLU & 6103.29 & 595.82 & 0.19 & 51.78 \\
\hline\end{tabular}\end{table}

\begin{table*}[!htbp]\centering
\caption{Pendulum results over both specified initial states. Blocks to tolerance are reported as $\xi_1/\xi_2$. The cost is from $\xi_1$. Times are mean/maximum in milliseconds.}
\label{tab:g2clock}
\begin{tabular}{llrrrrr}\hline
Algorithm&$\Vhat$ family&Blocks to tolerance&$J^{\rm ap}(NM)$&Inner solve&First command ready&Complete update\\\hline
3 & POLY & 89/83 & 234.20 & 16.82/37.93 & 16.82/37.94 & 16.82/37.94 \\
3 & PWA & 42/31 & 312.21 & 30.08/108.64 & 30.09/108.65 & 30.09/108.65 \\
3 & ReLU & 38/34 & 234.25 & 48.37/102.96 & 48.37/102.96 & 48.37/102.96 \\
4 & POLY & 89/83 & 234.35 & 27.01/66.83 & 27.06/66.86 & 27.09/66.90 \\
4 & PWA & 42/31 & 314.08 & 43.51/182.74 & 43.56/182.78 & 43.59/182.89 \\
4 & ReLU & 38/35 & 234.30 & 68.86/285.09 & 68.90/285.13 & 68.93/285.15 \\
\hline\end{tabular}\end{table*}

For $0\le k\le NM$, where $N$ is the number of complete blocks applied,
the accumulated executed cost on the measured trajectory is
\begin{equation}\label{eq:newexecutedcost}
J^{\rm ap}(k)=\sum_{s=0}^{k-1}\gamma^{\lfloor s/M\rfloor}
                 [V(x(s))+r(u(s))].
\end{equation}
Here $x(s)$ is the measured state before applying $u(s)$, and
$J^{\rm ap}(0)=0$. The exponent $\lfloor s/M\rfloor$ applies the
discount once per complete block as in~\eqref{eq:blockcost}.
The record ends at the first block boundary satisfying $|\theta|\le1.00$ mrad and $|\omega|\le0.01$ rad/s. 
Figures~\ref{fig:g2sets} to~\ref{fig:g2learning} and Table~\ref{tab:g2clock} report the 12 runs. All eight PWA and ReLU runs reach the tolerance within 200 samples. The four polynomial runs require 356 samples from $\xi_1$ and 332 from $\xi_2$. The replay checks all 2540 applied torques, 635 block acceptances, the fitted model used at each update, and the executed cost~\eqref{eq:newexecutedcost}. 
In every run the predicted balls contain the true successors and the state constraints hold. The reported costs accumulate up to the respective first tolerance hit.

Table~\ref{tab:g2clock} pools the measured updates from both initial states separately for each algorithm and value-function family. Inner solve includes numerical feasibility, objective minimization and verification of the returned torque vector. First command ready includes fitting and model acceptance for Algorithm~\ref{alg:modelupdate}.


\bibliographystyle{IEEEtran}
\bibliography{fsclf_model_based_rl_references}

\end{document}